\documentclass[final]{siamart171218}
\usepackage{amssymb}
\usepackage{amsmath}
\usepackage{mathdots}
\usepackage{microtype}
\usepackage{booktabs}
\usepackage{bm}
\usepackage{xcolor}
\usepackage[binary-units=true]{siunitx}
\usepackage[normalem]{ulem}
\usepackage{enumitem}
\usepackage{algorithm}
\usepackage{algorithmicx}
\usepackage[noend]{algpseudocode}
\makeatletter
\algrenewcommand\ALG@beginalgorithmic{\small}
\makeatother
\usepackage{tikz}
\usetikzlibrary{positioning}
\usetikzlibrary{shapes,arrows}
\usepackage{subcaption}
\usepackage{overpic}
\usepackage{import}

\renewcommand{\vec}[1]{\bm{#1}}
\newcommand{\trans}{\top}
\newcommand{\ldbrace}{\{\!\!\{}
\newcommand{\rdbrace}{\}\!\!\}}
\newcommand{\ldbrack}{[\![}
\newcommand{\rdbrack}{]\!]}
\newcommand{\avg}[1]{\ldbrace #1 \rdbrace}
\newcommand{\jump}[1]{\ldbrack #1 \rdbrack}
\newcommand{\flux}[1]{\widehat{#1}}
\protected\def\f{f}
\protected\def\c{c}
\DeclareMathOperator*{\argmin}{arg\,min}

\protected\def\squaresymbol{\scalebox{0.8}{\raisebox{0.08em}{\color[HTML]{0060ad}$\blacksquare$}}\hspace{0.2em}}
\protected\def\bulletsymbol{\scalebox{1.5}{\raisebox{-0.03em}{\color[HTML]{003866}$\bullet$}}\hspace{0.2em}}
\protected\def\trianglesymbol{\scalebox{1.1}[0.9]{\raisebox{0.1em}{\color[HTML]{008cff}$\blacktriangle$}}\hspace{0.2em}}
\protected\def\xsymbol{\scalebox{0.9}{\raisebox{0.08em}{$\bm{\times}$}}\hspace{0.2em}}
\protected\def\diamondsymbol{\raisebox{0.1em}{\scalebox{0.7}{\rotatebox[origin=c]{45}{\color[HTML]{4dafff}$\blacksquare$}}}\hspace{0.2em}}
\protected\def\triangledownsymbol{\scalebox{1.1}[0.9]{\raisebox{0.1em}{\color[HTML]{004680}$\blacktriangledown$}}\hspace{0.2em}}
\protected\def\plussymbol{\scalebox{0.9}{\raisebox{0.08em}{$\bm{+}$}}\hspace{0.2em}}
\protected\def\circlesymbol{\scalebox{0.7}{\raisebox{0.15em}{$\bigcirc$}}\hspace{0.2em}}

\numberwithin{equation}{section}
\numberwithin{theorem}{section}

\graphicspath{{figures/}}
\makeatletter
\def\input@path{{figures/}}
\makeatother

\title{Efficient Operator-Coarsening Multigrid Schemes for Local Discontinuous Galerkin Methods}
\author{Daniel Fortunato\thanks{Paulson School of Engineering and Applied Sciences, Harvard University, Cambridge, MA 02138 (\texttt{dfortunato@g.harvard.edu}, \texttt{chr@seas.harvard.edu}). This work is supported by the National Defense Science and Engineering Graduate Fellowship.} \and Chris H.~Rycroft\footnotemark[1] \and Robert Saye\thanks{Mathematics Group, Lawrence Berkeley National Laboratory, Berkeley, CA 94720 (\texttt{rsaye@lbl.gov})}. This research was supported by the Applied Mathematics Program of the U.S.~DOE Office of Advanced Scientific Computing Research under contract number DE-AC02-05CH11231. Some computations used resources of the National Energy Research Scientific Computing Center, a DOE Office of Science User Facility supported by the Office of Science of the U.S. Department of Energy under Contract No. DE-AC02-05CH11231.}
\date{\today}
\begin{document}
\maketitle

\begin{abstract}
An efficient $hp$-multigrid scheme is presented for local discontinuous Galerkin (LDG) discretizations of elliptic problems, formulated around the idea of separately coarsening the underlying discrete gradient and divergence operators. We show that traditional multigrid coarsening of the primal formulation leads to poor and suboptimal multigrid performance, whereas coarsening of the flux formulation leads to optimal convergence and is equivalent to a purely geometric multigrid method. The resulting operator-coarsening schemes do not require the entire mesh hierarchy to be explicitly built, thereby obviating the need to compute quadrature rules, lifting operators, and other mesh-related quantities on coarse meshes. We show that good multigrid convergence rates are achieved in a variety of numerical tests on 2D and 3D uniform and adaptive Cartesian grids, as well as for curved domains using implicitly defined meshes and for multi-phase elliptic interface problems with complex geometry. Extension to non-LDG discretizations is briefly discussed.
\end{abstract}



\section{Introduction}

Discontinuous Galerkin (DG) methods have gained broad popularity in recent years. They are well-suited to $hp$-adaptivity, provide high-order accuracy, and can be applied to a wide range of problems on complex geometries with unstructured meshes. Although DG methods were first applied to the discretization of hyperbolic conservation laws, they have been extended to handle elliptic problems and diffusive operators in a unified framework~\cite{Arnold_02_01}. Such methods include the symmetric interior penalty (SIP) method~\cite{Douglas_76_01, Arnold_82_01}, the Bassi--Rebay (BR1, BR2) methods~\cite{Bassi_97_01, Bassi_97_02}, the local discontinuous Galerkin (LDG) method~\cite{Cockburn_98_01}, the compact discontinuous Galerkin (CDG) method~\cite{Peraire_08_01}, the line-based discontinuous Galerkin method~\cite{Persson_13_01}, and the hybridizable discontinuous Galerkin (HDG) method~\cite{Cockburn_09_01}. In particular, the development of efficient solvers for DG discretizations of elliptic problems is an active area of research.

Among the panoply of DG methods for elliptic problems, the LDG method is a popular choice: it is accurate, stable, simple to implement, and extendable to higher-order derivatives~\cite{Yan_02_01}. Additionally, on Cartesian grids it has been shown to be superconvergent~\cite{Cockburn_01_01}. The LDG method results in symmetric positive (semi)definite discretizations which are well-suited to solution by efficient iterative methods. In particular, the multigrid method has emerged as a natural candidate due to its success in the continuous finite element and finite difference communities, both as a standalone solver and as a preconditioner for the conjugate gradient (PCG) method. However, direct application of standard multigrid techniques to DG discretizations of elliptic problems can result in suboptimal performance when inherited bilinear forms are employed~\cite{Antonietti_15_01, Gopalakrishnan_03_01}, and much work has gone into developing specialized smoothers and coarse-correction methods to remedy this issue, even on Cartesian grids~\cite{Kanschat_03_01, Fabien_17_01}.

When a mesh hierarchy is available, geometric $h$-multigrid is a natural choice of solver. Error estimates have been derived for a multilevel interior penalty (IP) method on unstructured meshes, yielding convergence factors in the range $\rho \approx 0.3$--$0.5$ for Poisson's equation~\cite{Gopalakrishnan_03_01, Brenner_05_01}; the method has also been applied to adaptively refined Cartesian grids with similar results~\cite{Kanschat_04_01}. Subsequent work describes how the multilevel IP method~\cite{Gopalakrishnan_03_01} can be used as a preconditioner for the LDG method---both for the Schur complement system (using conjugate gradient) and for the saddle-point system (using GMRES)---resulting in a bounded condition number with respect to mesh size $h$~\cite{Kanschat_03_01}. More recent work for LDG and IP uses a multigrid W-cycle on nested~\cite{Antonietti_15_01} and agglomerated~\cite{Antonietti_17_01} unstructured meshes; however, these results indicate poor convergence factors of $\rho \approx 0.8$--$0.9$ even with many smoothing steps. On non-nested polygonal meshes, $h$-independent iteration counts with convergence factors $\rho \approx 0.2$--$0.3$ are achieved for SIP using an additive Schwarz smoother with 3--8 smoothing steps per V-cycle~\cite{Antonietti_19_01}.

A popular choice for high-order DG methods is $p$- or $hp$-multigrid~\cite{Helenbrook_03_01, Fidkowski_05_01, Luo_06_01, Bassi_08_01}, where $p$ refers to coarsening the polynomial degree $p$ in the multilevel hierarchy, and $hp$ refers to some combined strategy of coarsening the mesh size $h$ as well as the polynomial degree $p$ of the underlying discretization. Using factor-of-two coarsening in $p$ with an element Jacobi smoother, convergence factors of $\rho \approx 0.5$ were achieved for Laplace's equation with $p \leq 4$~\cite{Helenbrook_03_01}. Fidkowski et al.~\cite{Fidkowski_05_01} used a line smoother with sequential coarsening in $p$ that gave similar results for convection--diffusion problems, though the performance degraded as $h \to 0$. A method employing an overlapping Schwarz smoother with factor-of-two coarsening was used with PCG on LDG discretizations up to $p=32$~\cite{Stiller_16_01}; this method exhibited good convergence factors of $\rho \leq 0.1$ on high-aspect-ratio Cartesian grids, at the cost of an expensive smoother.

Algebraic multigrid methods have also been applied to DG discretizations of elliptic problems. A hierarchy of operators can be defined by agglomerating neighboring unknowns based on smoothed aggregation (SA), resulting in average convergence factors of $\rho \approx 0.4$ and $\rho \approx 0.2$ for the bilinear BR2 and SIP methods, respectively~\cite{Prill_09_01}. An SA method employing energy minimization was used with PCG to achieve $h$-independent convergence factors of $\rho \approx 0.2$ for LDG discretizations, but performance degraded with increasing $p$~\cite{Olson_11_01}. A method based on unsmoothed aggregation was developed for the IP method using a coarse space consisting of continuous linear basis functions~\cite{Bastian_12_01}; this method proved robust for multi-phase problems with large jumps in ellipticity coefficient, but efficiency weakly degraded with mesh size. A related approach based on smoothed aggregation and low-order coarse grid correction yielded similar results~\cite{Siefert_14_01}.

Independent of the type of multigrid method, a particular fact to note---and something we believe underpins the difficulties in applying multigrid to DG---is that coarsening a fine-grid operator is not always the same as constructing that operator directly from the coarse grid. Indeed, it was noted by Antonietti et al.~\cite{Antonietti_15_01} that for all stable and strongly consistent DG methods, ``convergence cannot be independent of the number of levels if inherited bilinear forms are considered (i.e., the coarse solvers are the restriction of the stiffness matrix constructed on the finest grid).'' Furthermore they noted that non-inherited forms must be employed for the multigrid method to be scalable. In the context of two-level methods, the reason for this loss of scalability is known~\cite{Antonietti_12_01}. In this paper, for the LDG method, we present a simple modification to traditional multigrid operator coarsening that yields optimal multigrid convergence and can be extended to other DG discretizations of elliptic problems. We confirm that traditional coarsening of the fine-mesh elliptic operator results in poor performance, and show that the coarsening of the saddle-point flux formation restores optimal multigrid efficiency. Our approach is equivalent to pure geometric multigrid but avoids the need to explicitly build the coarse mesh and its associated components, such as quadrature rules, Jacobian mappings, lifting operators, and face-to-element enumerations---as discussed, this holds benefit for a variety of intricate DG implementations where building the coarse mesh can be problematic. Nevertheless we point out that in the pure geometric multigrid setting, quadrature-free DG methods~\cite{Antonietti_18_01} have recently been proposed which avoid the construction of coarse mesh quadrature rules.

The paper is structured as follows. In section~\ref{sec:DG}, we formulate a general DG discretization of Poisson's equation and derive the LDG method through the appropriate choice of numerical flux. In section~\ref{sec:multigrid}, we describe the construction of geometric $hp$-multigrid methods in the corresponding DG setting. In particular, we show that traditional operator coarsening can fail to create the coarse operator resulting from faithful rediscretization in a pure geometric multigrid setting for LDG methods, and present a modified coarsening strategy that remedies this. In section~\ref{sec:results}, we present numerical results for the standard and modified multigrid methods on uniform and adaptively-refined Cartesian grids in 2D and 3D. We conclude with some examples of multi-phase elliptic interface problems on implicitly defined meshes, which demonstrate good multigrid performance even on meshes with long and thin filaments as well as tiny and dispersed phase components.

\section{Discontinuous Galerkin formulation}\label{sec:DG}

\subsection{Model problem}
The model elliptic problem considered in this work is the Poisson problem
\begin{equation}\label{eq:poisson}
\begin{aligned}
-\nabla^2 u &= f &&\text{in } \Omega, \\
u &= g &&\text{on } \Gamma_D, \\
\nabla u \cdot \vec{n} &= h &&\text{on } \Gamma_N,
\end{aligned}
\end{equation}
where $\Omega$ is a domain in $\mathbb{R}^d$, $\Gamma_D$ and $\Gamma_N$ denote the components of $\partial \Omega$ on which Dirichlet and Neumann boundary conditions are imposed, $\vec{n}$ is the outward unit normal to the boundary, and $f$, $g$, and $h$ are given functions defined on $\Omega$ and its boundary.

\subsection{DG for elliptic problems}
In order to apply a DG method to \eqref{eq:poisson}, we rewrite it as a first-order system by introducing the auxiliary variable $\vec{q} = \nabla u$ and writing the Laplacian as the divergence of $\vec{q}$~\cite{Arnold_02_01}:
\begin{equation}\label{eq:poisson_fo}
\begin{aligned}
\vec{q} &= \nabla u &&\text{in } \Omega, \\
-\nabla \cdot \vec{q} &= f &&\text{in } \Omega, \\
u &= g &&\text{on } \Gamma_D, \\
\vec{q} \cdot \vec{n} &= h &&\text{on } \Gamma_N.
\end{aligned}
\end{equation}
In this work, we mainly consider discretizations of \eqref{eq:poisson_fo} wherein the corresponding meshes arise from Cartesian grids, quad/octrees, or implicitly defined meshes of more complex curved domains (see sections \ref{sec:cartesian}, \ref{sec:amr}, and \ref{sec:implicit_mesh}, respectively). As such, it is natural to adopt a tensor-product piecewise polynomial space. Let $\mathcal E = \bigcup_i E_i$ denote the set of elements of a mesh of $\Omega$, let $p \geq 1$ be an integer, and define $\mathcal Q_p(E)$ to be the space of tensor-product polynomials of degree $p$ on the element $E$. For example, $\mathcal Q_3$ is the space of bicubic (in 2D) or tricubic (in 3D) polynomials having 16 or 64 degrees of freedom, respectively.
We define the corresponding spaces of discontinuous piecewise polynomials and vector fields on the mesh as
\begin{align}
V_h(\mathcal{E}) &= \big\{ v : \Omega \to \mathbb{R} \enskip\bigl|\enskip v|_E \in {\mathcal Q}_p(E) \text{ for every } E \in {\mathcal E} \bigr\}, \\
V_h^d(\mathcal{E}) &= \big\{ \vec{\omega} : \Omega \to \mathbb{R}^d \enskip\bigl|\enskip \vec{\omega}|_E \in [{\mathcal Q}_p(E)]^d \text{ for every } E \in {\mathcal E} \bigr\},
\end{align}
respectively. We denote by $(\cdot,\cdot)$ the natural $L^2$ inner product on $V_h$ and by $\| \cdot \|$ the corresponding norm, $\|u\|^2 = (u,u)$, with analogous definitions for $V_h^d$.

In a DG method, both $\vec{q}$ and its divergence are defined weakly via numerical fluxes defined on each mesh face. The weak form of \eqref{eq:poisson_fo} consists of finding $(\vec{q}_h,u_h) \in V_h^d \times V_h$ such that
\begin{gather}
\int_E \vec{q}_h \cdot \vec{\omega} = -\int_E u_h \nabla \cdot \vec{\omega} + \int_{\partial E} \flux{u}_h\,\vec{\omega} \cdot \vec{n}, \label{eq:grad_weak1} \\[0.4em]
\int_E \vec{q}_h \cdot \nabla v - \int_{\partial E} \flux{\vec{q}}_h\,v \cdot \vec{n} = \int_E f\,v, \label{eq:div_weak1}
\end{gather}
for all test functions $(\vec{\omega},v) \in [\mathcal{Q}_p(E)]^d \times \mathcal{Q}_p(E)$ and for all $E \in \mathcal{E}$. The numerical fluxes $\flux{\vec{q}}_h$ and $\flux{u}_h$ are approximations to $\vec{q}_h$ and $u_h$, respectively, on each mesh face and define how the degrees of freedom in each element are coupled together.

To more succinctly describe the coupling between elements, the following standard notation is adopted. Consider two adjacent elements $E^+$ and $E^-$ which share a face in $\mathcal{E}$. Let $\vec{n}^\pm$ denote the outward unit normals of $\partial E^\pm$ along the shared face and $(\vec{\omega}^\pm, v^\pm)$ denote the traces of $(\vec{\omega},v) \in V_h^d \times V_h$ from $E^\pm$ on the shared face. The average $\avg{\cdot}$ and jump $\jump{\cdot}$ operators on the shared face are then defined as
\begin{alignat*}{4}
\avg{&\vec{\omega}} &&= \tfrac12 (\vec{\omega}^+ + \vec{\omega}^-), &\qquad \avg{&v} &&= \tfrac12 (v^+ + v^-), \\
\jump{&\vec{\omega}} &&= \vec{\omega}^+ \cdot \vec{n}^+ + \vec{\omega}^- \cdot \vec{n}^-, &\qquad \jump{&v} &&= v^+ \vec{n}^+ + v^- \vec{n}^-.
\end{alignat*}
On boundary faces, $(\vec{\omega}^-, v^-)$ shall refer to the traces of $(\vec{\omega},v)$ from the corresponding element touching $\partial \Omega$.

\subsection{The local discontinuous Galerkin method}\label{sec:LDG}
The choice of numerical flux in \eqref{eq:grad_weak1}--\eqref{eq:div_weak1} defines a DG method. Here we focus on the LDG method~\cite{Cockburn_98_01}, which chooses numerical fluxes $\flux{\vec{q}}_h$ and $\flux{u}_h$ according to the general form
\begin{equation}\label{eq:q_flux}
\flux{\vec{q}}_h = \left\{
\begin{array}{cl}
\avg{\vec{q}_h} + \vec{\beta} \jump{\vec{q}_h} - \tau_0\jump{u_h} & \text{on any interior face,} \\
\vec{q}_h^- -\tau_D (u_h^- - g) \vec{n} & \text{on any face of $\Gamma_D$,} \\
h\,\vec{n} & \text{on any face of $\Gamma_N$},
\end{array}
\right.
\end{equation}
and
\begin{equation}\label{eq:u_flux}
\flux{u}_h = \left\{
\begin{array}{cl}
\avg{u_h} - \vec{\beta}\cdot\jump{u_h} & \text{on any interior face,} \\
g & \text{on any face of $\Gamma_D$,} \\
u_h^- & \text{on any face of $\Gamma_N$},
\end{array}
\right.
\end{equation}
where $\vec{\beta}$ is a (possibly face-dependent) user-defined vector; for example, in a one-sided flux scheme, $\vec{\beta} = \pm \tfrac12 \vec{n}$. Here, the numerical flux $\flux{\vec{q}}_h$ includes penalty stabilization terms; $\tau_0 \geq 0$ is a penalty parameter associated with interior faces and $\tau_D > 0$ is associated with Dirichlet boundary faces (if any). Generally, $\tau_0$ must be strictly positive to ensure well-posedness of the discrete problem, but in some cases (e.g., on Cartesian grids with particular choices of $\vec{\beta}$), $\tau_0$ can be set equal to zero~\cite{Cockburn_07_01}. If $\Gamma_D$ is nonempty, $\tau_D$ must be positive to ensure well-posedness of the final discrete problem. To be consistent with the scaling of penalty parameters in other DG methods, we choose the penalty parameters to scale inversely with the element size $h$~\cite{Castillo_00_01} so that $\tau_0 = \widetilde{\tau}_0/h$ and $\tau_D = \widetilde{\tau}_D/h$ where $\widetilde{\tau}_0$ and $\widetilde{\tau}_D$ are constants.\footnote{To obtain uniform stability in the limit of large $p$, $\tau_0$ should also scale with $p^2$~\cite{Perugia_02_01}; we do not consider this aspect for the moderate values of $p$ tested in this work ($p=1\mbox{--}8$).} Furthermore, although arbitrarily small choices of $\widetilde{\tau}_0$ and $\widetilde{\tau}_D$ suffice to ensure well-posedness, later we show that a carefully considered choice of these values can greatly benefit multigrid performance (see section~\ref{sec:tau_study}).

We first particularize \eqref{eq:grad_weak1} for the LDG method, which is the weak statement that $\vec{q} = \nabla u$. We slightly modify the weak form \eqref{eq:grad_weak1} by defining $\vec{q}_h \in V_h^d$ in strong-weak form,\footnote{The \textit{strong-weak form} states that $\vec{q}_h$ must satisfy $\int_E \vec{q}_h \cdot \vec{\omega} = \int_E \nabla u_h \cdot \vec{\omega} + \int_{\partial E} (\flux{u}_h - u_h)\, \vec{\omega} \cdot \vec{n}$ whereas the \textit{weak form} states that $\vec{q}_h$ must satisfy $\int_E \vec{q}_h \cdot \vec{\omega} = -\int_E u_h\,\nabla \cdot \vec{\omega} + \int_{\partial E} \flux{u}_h \, \vec{\omega} \cdot \vec{n}$. The two forms are equivalent whenever the employed quadrature scheme exactly satisfies the identity of integration by parts, which in practice is generally true for quadrilateral, prismatic, simplicial elements, etc., but is generally not true when approximate numerical quadrature schemes are used, e.g., as on implicitly defined curved elements. In the latter situation, to ensure symmetry of the final discrete Laplacian operator, it is necessary to use the strong-weak form to define $\vec{q}_h$ and the weak form to define the divergence of $\vec{q}_h$ (or vice versa)~\cite{Saye_17_01}.} such that
\begin{equation}\label{eq:qdef}
\int_E \vec{q}_h \cdot \vec{\omega} = \int_E \nabla u_h \cdot \vec{\omega} + \int_{\partial E} (\flux{u}_h - u_h)\, \vec{\omega} \cdot \vec{n}
\end{equation}
holds for every element $E \in {\mathcal E}$ and every test function $\vec{\omega} \in [\mathcal{Q}_p(E)]^d$. Upon summing \eqref{eq:qdef} over every element of the mesh and using the definition of the numerical flux $\flux{u}_h$ in \eqref{eq:u_flux}, we have that, for any $\vec \omega \in V_h^d$,
\begin{equation}\label{eq:qdefintegral}
\int_\Omega \vec{q}_h \cdot \vec{\omega} = \sum_{E \in \mathcal E} \int_E \nabla u_h \cdot \vec{\omega} - \int_{\Gamma_0} \jump{u_h} \cdot \left(\avg{\vec{\omega}} + \vec{\beta} \jump{\vec{\omega}}\right) + \int_{\Gamma_D} (g - u_h^-) \vec{\omega}^- \cdot \vec{n},
\end{equation}
where $\Gamma_0$ denotes the union of all interior faces of $\mathcal{E}$.
Define the following operators:
\begin{itemize}[itemsep=0.5em]
\item Let $\nabla_h : V_h \to V_h^d$ be the \emph{broken gradient operator} and $L : V_h \to V_h^d$ be the \emph{lifting operator}, such that
\begin{align*}
&\int_\Omega (\nabla_h u) \cdot \vec{\omega} = \sum_{E \in \mathcal E} \int_E \nabla u \cdot \vec \omega, \\
&\int_\Omega (Lu) \cdot \vec{\omega} = -\int_{\Gamma_0} \jump{u} \cdot \left(\avg{\vec{\omega}} + \vec{\beta} \jump{\vec{\omega}}\right) - \int_{\Gamma_D} u^- \vec{\omega}^- \cdot \vec{n}
\end{align*}
holds for every $\vec \omega \in V_h^d$ and each $u \in V_h$.
\item Define $J_D(g) \in V_h^d$ such that
\[
\int_\Omega J_D(g) \cdot \vec{\omega} = \int_{\Gamma_D} g\,\vec{\omega}^- \cdot \vec{n}
\]
holds for every $\vec{\omega} \in V_h^d$.
\end{itemize}
Accordingly, \eqref{eq:qdefintegral} is equivalent to the statement that
\begin{equation}\label{eq:qresult}
\vec{q}_h = (\nabla_h + L)u_h + J_D(g) = G u_h + J_D(g)
\end{equation}
where $G : V_h \to V_h^d$ is the \emph{discrete gradient operator}, $G = \nabla_h + L$. The formula \eqref{eq:qresult} is the LDG discretization of the statement $\vec q = \nabla u$, taking into account Dirichlet boundary data.

Next, we particularize \eqref{eq:div_weak1} for the LDG method, which is the weak statement that $-\nabla \cdot \vec{q} = f$. Upon summing \eqref{eq:div_weak1} over every mesh element and using the definition of the numerical flux $\flux{\vec{q}}_h$ in \eqref{eq:q_flux}, we have that, for any $v \in V_h$,
\begin{equation}\label{eq:wdef2}
\begin{aligned}
\sum_{E \in \mathcal E} \int_E \vec{q}_h \cdot \nabla v &- \int_{\Gamma_0} (\avg{\vec{q}_h} + \vec{\beta} \jump{\vec{q}_h} - \tau_0 \jump{u_h}) \cdot \jump{v} - \int_{\Gamma_D} (\vec{q}_h^- \cdot \vec{n} - \tau_D u_h^-)\,v^- \\
&= \int_\Omega f\,v + \tau_D \int_{\Gamma_D} g\,v^- + \int_{\Gamma_N} h\,v^-.
\end{aligned}
\end{equation}
Additionally, define the following operators:
\begin{itemize}[itemsep=0.5em]
\item Similar to the operator $J_D$ above, let $J_N(h) \in V_h$ be such that
\[
\int_\Omega J_N(h)\,v = \int_{\Gamma_N} h\,v^-
\]
for all $v \in V_h$.
\item Let $E_0, E_D : V_h \to V_h$ be the operators such that, for each $u \in V_h$,
\[
\int_\Omega E_0 (u)\, v = \int_{\Gamma_0} \jump{u} \cdot \jump{v}, \qquad \int_\Omega E_D (u)\,v = \int_{\Gamma_D} u^- v^-
\]
hold for every $v \in V_h$. These operators penalize jumps in the discrete solution on interior and Dirichlet boundary faces, respectively.
\item Let $a_D(g) \in V_h$ be such that
\[
\int_{\Omega} a_D(g)\, v = \int_{\Gamma_D} g\,v^-
\]
for all $v \in V_h$.
\end{itemize}
Then, using the fact that $(\vec{q}_h, \nabla_h v) + (\vec{q}_h, Lv) = (\vec{q}_h, Gv)$, \eqref{eq:wdef2} is equivalent to 
\begin{equation}\label{eq:wresult}
(\vec{q}_h, Gv) + \tau_0 (E_0 u_h, v) + \tau_D (E_D u_h, v) = (f + J_N(h) + \tau_D a_D(g), v),
\end{equation}
or, putting aside penalty terms, $G^\ast \vec{q}_h = f + J_N(h)$, where $-G^\ast$ is the \emph{discrete divergence operator}, the negative adjoint of the discrete gradient operator $G$; this is the LDG discretization of the statement that $-\nabla \cdot \vec{q} = f$, taking into account Neumann boundary data.

\subsubsection{Primal formulation}
To obtain the \emph{primal formulation} of the LDG method, we combine \eqref{eq:wresult} with \eqref{eq:qresult} to eliminate $\vec{q}_h$ and arrive at an equation for $u_h$. The primal LDG formulation of \eqref{eq:poisson} reads as follows: find $u_h \in V_h$ such that
\begin{equation}\label{eq:primal}
a(u_h,v) = \ell(v)
\end{equation}
for all $v \in V_h$, where the bilinear form $a(\cdot,\cdot)$ is given by
\[
a(u_h,v) = (Gu_h, Gv) + \tau_0 (E_0 u_h, v) + \tau_D (E_D u_h, v)
\]
and the linear functional $\ell(\cdot)$ is given by
\[
\ell(v) = (f, v) - (J_D(g), Gv) + (J_N(h), v) + \tau_D (a_D(g), v).
\]
One may verify that the bilinear form $a(u,v)$ is symmetric. Discretization of the primal form \eqref{eq:primal} with respect to a particular basis of $V_h$ yields a symmetric positive (semi)definite linear system of the form\footnote{Throughout the paper we shall frequently use the same symbol to denote (i) elements of spaces such as $V_h$ or operators acting on such elements, and (ii) vectors of coefficients in the chosen basis or matrices acting on such vectors. The distinction should be clear from context. Further comments are provided in section \ref{sec:choicebasis}.}
\begin{equation}\label{eq:primal_disc}
A u_h = \ell,
\end{equation}
where $A$ is the matrix form of the negative \emph{discrete Laplacian operator}.

\subsubsection{Flux formulation}
An alternative, but equivalent, characterization of the LDG method is the so-called \emph{flux formulation}, which does not eliminate the auxiliary variable $\vec{q}_h$ from the system \eqref{eq:qresult} and \eqref{eq:wresult} but instead retains it as a primary unknown. The flux formulation of \eqref{eq:poisson} then reads as follows: find $(\vec{q}_h,u_h) \in V_h^d \times V_h$ such that
\begin{equation}\label{eq:flux}
\begin{aligned}
m(\vec{q}_h,\vec{\omega}) - \operatorname{grad}(u_h,\vec{\omega}) &= j(\vec{\omega}), \\
-\operatorname{div}(\vec{q}_h,v) + \tau(u_h,v) &= k(v),
\end{aligned}
\end{equation}
for all $(\vec{\omega}, v) \in V_h^d \times V_h$, where
\begin{align*}
m(\vec{q},\vec{\omega}) &= (\vec{q},\vec{\omega}), \\
\operatorname{grad}(u,\vec{\omega}) &= (Gu,\vec{\omega}), \\
\operatorname{div}(\vec{q},v) &= -(\vec{q},Gv), \\
\tau(u,v) &= \tau_0 (E_0 u,v) + \tau_D (E_D u, v), \\
j(\vec{\omega}) &= (J_D(g),\vec{\omega}), \\
k(v) &= (f,v) + (J_N(h),v) + \tau_D (a_D(g),v).
\end{align*}
Discretization of the flux form \eqref{eq:flux} with respect to a particular basis of $V_h^d \times V_h$ yields a symmetric positive (semi)definite linear system of the form
\begin{equation}\label{eq:LDG_block}
\begin{bmatrix}
 M   & -M G \\[.2em]
-M D & \phantom{-}M T
\end{bmatrix}
\begin{bmatrix}
\vec{q}_h \\[.2em] u_h
\end{bmatrix}
=
\begin{bmatrix}
j \\[.2em] k
\end{bmatrix}.
\end{equation}
Here, $M$ is the block diagonal mass matrix for $V_h$ or $V_h^d$ (depending on context), $G$ is the matrix form of the discrete gradient operator, $D = -M^{-1} G^\trans M$ is the matrix form of the discrete divergence operator, and $T = \tau_0 E_0 + \tau_D E_D$ contains the discrete penalty terms. Since $M^{-1}$ is also block diagonal, we can easily take the Schur complement of $M$ in \eqref{eq:LDG_block} to obtain a linear system for the unknown vector $u_h$,
\begin{equation}\label{eq:LDG_schur}
A u_h = \ell,
\end{equation}
where $A = M (-D G + T) = G^\trans M G + M T$ and $\ell = k - G^\trans j$.

The reduced linear system \eqref{eq:LDG_schur} is equivalent to the discrete primal formulation \eqref{eq:primal_disc}. However, as we will demonstrate next in section~\ref{sec:multigrid}, the two formulations have different implications for multigrid methods. In particular, applying standard operator coarsening to the discrete primal formulation results in poor multigrid performance; coarsening the discrete flux formulation \eqref{eq:LDG_block} in both $\vec{q}_h$ and $u_h$ \emph{before} taking the Schur complement \eqref{eq:LDG_schur} results in optimal multigrid performance, and is equivalent to pure geometric multigrid.

\subsubsection{Remarks on the choice of basis}\label{sec:choicebasis}

The analysis and discussion presented in this paper is agnostic to the particular choice of basis for the piecewise polynomial space $V_h$. One may use a nodal basis, a modal basis, or some other choice, provided it is understood that every basis-dependent matrix (e.g., the mass matrix $M$) is defined consistently, relative to the chosen basis. In a numerical implementation, one should consider aspects of conditioning, accuracy, stability, sparsity, and computational complexity. For example, for low-to-moderate polynomial degree on rectangular elements, as used in this work, a tensor-product nodal basis using Gauss--Lobatto nodes is a natural choice~\cite{Hesthaven_08_01}; for very large $p$, a modal basis may have better conditioning or improved cost of mass matrix inversion, and thus may be more suitable. In our particular implementation, we have used a tensor-product Gauss--Lobatto nodal basis. We emphasize however that the presented multigrid methods and the essential conclusions drawn are not dependent on this choice.

\section{Multigrid methods}\label{sec:multigrid}
We assume here that the reader has some familiarity with multigrid methods; see for example Briggs, Henson, and McCormick~\cite{Briggs_00_01} for a review of their design and operation. A geometric multigrid method consists of four main ingredients: a mesh hierarchy, an interpolation operator to transfer approximate solutions from a coarse mesh onto a fine mesh, a restriction operator to formulate a coarse mesh correction problem by restricting the residual from the fine mesh, and a smoother/relaxation method. We consider these ingredients separately first, and then combine them into a multigrid V-cycle.

Multigrid methods rely on the complementarity between relaxation and interpolation. In the geometric multigrid context, a relaxation method that is effective at damping high-frequency, oscillatory errors but slow to damp smooth, low-frequency ones benefits from the action of an interpolation operator that can accurately transfer low-frequency information. By solving a correction equation for the error on a coarser grid, fine-grid low-frequency errors become coarse-grid high-frequency errors for which coarse-grid relaxation is effective. An interpolation operator then transfers this low-frequency correction to the fine grid.

In the following sections we focus our description on $h$-multigrid methods, wherein the mesh is coarsened geometrically at each level. However, much of our analysis carries over to $p$-multigrid methods, which hold the mesh fixed and instead coarsen the polynomial space by reducing $p$ at each level. We will try to point out the distinctions between the two methods when the analogues are not immediately obvious, though we will use the notation $h$ in our description.

\subsection{Mesh hierarchy}
In this work, we employ quadtrees (in 2D) and octrees (in 3D) to define the finest mesh---whether it is uniform, adaptively refined, or used as the background grid for an implicitly defined mesh (see section \ref{sec:implicit_mesh}). The tree structure then naturally defines a hierarchy of nested meshes for use in $h$-multigrid, which are spatially coarsened by a factor of two in each dimension on each level. For adaptively refined meshes where the cell size is not uniform, we coarsen each element as rapidly as the tree structure permits (see, e.g., Figure~\ref{fig:amr}).

In the context of $p$-multigrid methods, a mesh hierarchy is defined by applying a specific $p$-coarsening strategy to the fine mesh. For example, one could coarsen $p$ sequentially ($p \to p-1 \to p-2 \to \ldots \to 1$), by a factor of two ($p \to p/2 \to p/4 \to \cdots \to 1$), or by some user-defined sequence of $p$'s. The first method is a common choice when low-order polynomials are used on the finest mesh, whereas the second method is better suited to high-order discretizations.

Mesh hierarchies can also be generated by combining coarsening strategies in both $h$ and $p$. For example, a popular choice is to layer $p$-multigrid on top of $h$-multigrid, so that $h$-multigrid with a low-order polynomial degree is used as the bottom solver in the $p$-multigrid hierarchy.

\subsection{Interpolation}
The interpolation operator $I_\c^\f$ transfers a piecewise polynomial function $u_\c \in V_{2h}(\mathcal{E}_\c)$ defined on a coarse mesh $\mathcal{E}_\c$ to a piecewise polynomial function $u_\f \in V_h(\mathcal{E}_\f)$ on a fine mesh $\mathcal{E}_\f$. (Throughout this work, subscripts or superscripts $\f$, $h$ and $\c$, $2h$ shall denote objects corresponding to the fine mesh and coarse mesh, respectively.) We define the interpolation operator so that it injects the piecewise polynomial function on the coarse mesh into the fine mesh, unmodified. From the $h$-multigrid perspective, $u_\f|_{E_\f}$ on the fine element $E_\f$ is simply the polynomial $u_\c|_{E_\c}$ restricted to $E_\f$, where $E_\c \supset E_\f$ is the corresponding coarse element in the mesh hierarchy. From the $p$-multigrid perspective, the lower-degree polynomial $u_\c$ can be exactly represented as a higher-degree polynomial by taking the higher-order coefficients of $u_\f$ to be zero. In either case, when regarded as an operator from \smash{$L^2(\Omega) \to L^2(\Omega)$}, \smash{$I_\c^\f$} is the identity operator. The operator is linear and has the property that it preserves constant functions, i.e., $u_\c \equiv 1$ is mapped to \smash{$u_\f \equiv 1$}. This property ensures that, throughout a V-cycle, the coarse mesh discrete problems preserve the compatibility condition required in semidefinite problems having solely Neumann boundary conditions.

\subsection{Restriction}
We define the restriction operator $R_\f^\c : V_h(\mathcal{E}_\f) \to V_{2h}(\mathcal{E}_\c)$ as the adjoint of the interpolation operator, i.e., such that
\begin{equation}\label{eq:rdef}
(R_\f^\c u_\f, u_\c)_{\mathcal{E}_\c} = (u_\f, I_\c^\f u_\c)_{\mathcal{E}_\f}
\end{equation}
holds for every $u_\f \in V_h(\mathcal{E}_\f)$ and every $u_\c \in V_{2h}(\mathcal{E}_\c)$. Equivalently, letting $I_\c^\f$, $R_\f^\c$, $u_\c$, and $u_\f$ also denote matrices and vectors relative to the user-defined bases of $V_h(\mathcal{E}_\f)$ and $V_{2h}(\mathcal{E}_\c)$, \eqref{eq:rdef} can be restated as
\[
(R_\f^\c u_\f)^\trans M_\c u_\c = u_\f^\trans M_\f I_\c^\f u_\c
\]
where $M_\c$ and $M_\f$ are the block-diagonal mass matrices of the coarse and fine meshes, respectively. Therefore,
\begin{equation}\label{eq:rmatrix}
R_\f^\c = M_\c^{-1} (I_\c^\f)^\trans M_\f.
\end{equation}
Defining the restriction operator in this manner---sometimes referred to as Galerkin projection---results in several notable properties:
\begin{itemize}[itemsep=0.5em]
\item One may interpret $R_\f^\c u_\f$ as ``averaging'' elemental polynomials of $u_\f \in V_h(\mathcal{E}_\f)$ on the fine mesh to determine a coarsened piecewise-polynomial representation on the coarse mesh. The averaging is performed in a way that locally preserves the mass of $u_\f$: indeed, since $I_\c^\f$ preserves constant functions, we have that $(R_\f^\c u_\f, 1)_{\mathcal{E}_\c} = (u_\f, 1)_{\mathcal{E}_\f}$ for all $u_\f \in V_h(\mathcal{E}_\f)$.
\item The adjoint method can also be viewed as an $L^2$ projection of $u_\f \in V_h(\mathcal{E}_\f)$ onto $V_{2h}(\mathcal{E}_\c)$. The variational problem $\argmin_{u_\c \in V_{2h}(\mathcal{E}_\c)} \|u_\c - u_\f\|_{\Omega}^2$ optimizes the functional
\[
V_{2h}(\mathcal{E}_\c) \ni u_\c \mapsto (u_\c,u_\c)_{\mathcal{E}_\c} - 2 (u_\f, I_\c^\f u_\c)_{\mathcal{E}_\f} = u_\c^\trans M_\c u_\c - 2 u_\f^\trans M_\f I_\c^\f u_\c
\]
whose unique minimum is given by $u_\c = M_\c^{-1} (I_\c^\f)^\trans M_\f u_\f$.
\item From the preceding property, one can immediately infer that
\begin{equation}\label{eq:riid}
R_\f^\c I_\c^\f = \mathbb{I},
\end{equation}
where $\mathbb{I}$ is the identity operator; i.e., interpolating a piecewise polynomial function from a coarse mesh onto a fine mesh and immediately restricting the result shall return the original function. The relation in \eqref{eq:riid} together with \eqref{eq:rmatrix} also provides a method to compute the coarse-mesh mass matrix from the fine-mesh mass matrix:
\begin{equation}\label{eq:coarse_mass}
M_\c = (I_\c^\f)^\trans M_\f I_\c^\f.
\end{equation}
\end{itemize}

Given a linear operator $A : V_h(\mathcal{E}_\f) \to V_h(\mathcal{E}_\f)$, one can define a coarsened operator $\mathcal{C}(A) : V_{2h}(\mathcal{E}_\c) \to V_{2h}(\mathcal{E}_\c)$ in a similar way, by proceeding variationally: we define $\mathcal{C}(A)$ such that
\[
(\mathcal{C}(A) u_\c, v_\c)_{\mathcal{E}_\c} = (A I_\c^\f u_\c, I_\c^\f v_\c)_{\mathcal{E}_\f}
\]
holds for all $u_\c, v_\c \in V_{2h}(\mathcal{E}_\c)$. Viewing $A$ and $\mathcal{C}(A)$ as matrix operators, mapping vectors in the user-defined bases of $V_h(\mathcal{E}_\f)$ and $V_{2h}(\mathcal{E}_\c)$,
\[
\mathcal{C}(A) = M_\c^{-1} (I_\c^\f)^\trans M_\f A I_\c^\f = R_\f^\c A I_\c^\f.
\]
The last form is perhaps more commonly seen or referred to as \textit{``RAT''} in the multigrid literature~\cite{Xu_02_01}, where \textit{R} is restriction, \textit{A} is the fine-mesh operator, and \textit{T} (or \textit{P}) is the interpolation (or prolongation) operator; the essence of the present work is to show that directly applying \textit{RAT} to the negative discrete Laplacian resulting from the primal formulation of an LDG method results in an inefficient multigrid algorithm and that, instead, applying \textit{RAT} to the flux formulation, $\vec q = \nabla u$, $-\nabla \cdot \vec q = f$, leads to more efficient multigrid solvers.

\subsection{Operator coarsening and pure geometric multigrid}
In this section we compare a standard, purely geometric multigrid method to two multigrid schemes based on operator coarsening: (i) applying \textit{RAT} to the negative discrete Laplacian matrix of the primal formulation (``primal coarsening'') and (ii) applying \textit{RAT} to the $2\times2$ block matrix of the flux formulation (``flux coarsening''). By a pure geometric method, we mean one in which each level of the hierarchy is explicitly meshed and the LDG formulation is canonically applied to each level, with the above restriction and interpolation operators transferring residual and correction vectors (in a V-cycle) between levels. Our motivation here concerns an $h$-multigrid method; however, much of the following discussion has direct analogy with $p$-multigrid methods. In addition, in the context of DG methods requiring penalty parameters, a design choice can be made as to how the value of the penalty parameter is chosen at each level of the hierarchy. In this work we consider the natural choice in which every level of the hierarchy inherits the same value as the finest mesh. With this in mind, we discuss interaction between a pair of levels: suppose $\mathcal{E}_\f$ is the mesh of a fine level and $\mathcal{E}_\c$ is the mesh of the next-coarsest level.

\subsubsection{Primal coarsening}
Recall the primal form of the negative discrete Laplacian operator of an LDG method: as a matrix mapping the coefficient vectors in the basis of $V_h(\mathcal{E}_\f)$ into the basis of $V_h(\mathcal{E}_\f)$, i.e., premultiplying \eqref{eq:LDG_schur} by $M^{-1}$,
\[
A = -D G + \tau_0 E_0 + \tau_D E_D,
\]
where $G = \nabla_h + L$ is the discrete gradient operator and $D = -M^{-1} G^\trans M$ is the discrete divergence operator. To discuss the application of \textit{RAT} to $A$ and how it relates to a geometric multigrid implementation, we consider the individual terms making up $A$.

\begin{itemize}[itemsep=0.5em]
\item First, we note that the broken gradient operator satisfies the \textit{RAT} property, i.e., $\mathcal{C}(\nabla_h) = \nabla_{2h}$. Computing the piecewise gradient on a coarse mesh and interpolating the result to the fine mesh is the same as computing the piecewise gradient of the interpolant, i.e., $I_\c^\f \nabla_{2h} u_\c = \nabla_h I_\c^\f u_\c$ for all $u_\c \in V_{2h}(\mathcal{E}_\c)$; consequently, $\mathcal{C}(\nabla_h) = R_\f^\c \nabla_h I_\c^\f = R_\f^\c I_\c^\f \nabla_{2h} = \nabla_{2h}$ by \eqref{eq:riid}.
\item The lifting operator also satisfies the \textit{RAT} property, i.e., $\mathcal{C}(L_\f) = L_\c$. This is perhaps not immediately obvious, since source terms on a coarse mesh face will lift into the corresponding large coarse element, whereas the corresponding source terms on the fine mesh faces lift only into the smaller elements touching that face; however, the restriction of the result on the set of smaller elements agrees with the result of $L_\c$. To see this, apply the variational formulation of $\mathcal{C}(\cdot)$ to observe that
\begin{align*}
& (\mathcal{C}(L_\f) u_\c, \vec{v}_\c)_{\mathcal{E}_\c} \\
&~~= (L_\f I_\c^\f u_\c, I_\c^\f \vec{v}_\c)_{\mathcal{E}_\f} \\
&~~= -\int_{\Gamma_{0,\f}} \jump{I_\c^\f u_\c} \cdot \left( \avg{I_\c^\f \vec{v}_\c} + \vec{\beta} \jump{I_\c^\f \vec{v}_\c} \right) - \int_{\Gamma_{D,\f}} (I_\c^\f u_\c)^- (I_\c^\f \vec{v}_\c)^- \cdot \vec{n} \\
&~~= -\int_{\Gamma_{0,\c}} \jump{u_\c} \cdot \left( \avg{\vec{v}_\c} + \vec{\beta} \jump{\vec{v}_\c} \right) - \int_{\Gamma_{D,\c}} u_\c^- \vec{v}_\c^- \cdot \vec{n} \\
&~~= (L_\c u_\c, \vec{v}_\c)_{\mathcal{E}_\c}
\end{align*}
holds for all $u_\c \in V_{2h}(\mathcal{E}_\c)$ and $\vec{v}_\c \in V_{2h}^d(\mathcal{E}_\c)$. Here, $\Gamma_{0,\f}$ and $\Gamma_{0,\c}$ denote the union of interior faces of the fine and coarse meshes, respectively, and similarly for $\Gamma_{D,\f}$ and $\Gamma_{D,\c}$. The third equality holds because the interpolation operator introduces no nonzero jumps on the set of new fine mesh faces, i.e., on $\Gamma_{0,\f} \setminus \Gamma_{0,\c}$ and $\Gamma_{D,\f} \setminus \Gamma_{D,\c}$. (The preceding assumes that fine mesh faces inherit the same $\vec{\beta}$ value as coarse mesh faces; in particular, this is true for the one-sided LDG scheme in which $\vec{\beta} = \pm \tfrac12 \vec{n}$.)
\item It immediately follows from the preceding two properties that $\mathcal{C}(G_\f) = G_\c$. Moreover, 
\begin{align*}
\mathcal{C}(D_\f) &= R_\f^\c D_\f I_\c^\f = \bigl(M_\c^{-1} (I_\c^\f)^\trans M_\f \bigr) \bigl(-M_\f^{-1} G_\f^\trans M_\f \bigr) I_\c^\f \\
&= -M_\c^{-1} \bigl((I_\c^\f)^\trans G_\f^\trans M_\f I_\c^\f M_\c^{-1} \bigr) M_\c = -M_\c^{-1} \bigl(\mathcal{C}(G_\f) \bigr)^\trans M_\c \\
&= -M_\c^{-1} G_\c^\trans M_\c = D_\c.
 \end{align*}
\item It is straightforward to show that the penalty operators also satisfy the \textit{RAT} property, i.e., $\mathcal{C}(E_{0,\f}) = E_{0,\c}$ and $\mathcal{C}(E_{D,\f}) = E_{D,\c}$. As in the case of the lifting operator, this property derives from the fact the interpolation operator does not introduce jumps on fine mesh faces that do not overlap with coarse mesh faces.
\end{itemize}

Despite these consistencies, the negative discrete Laplacian does not satisfy the \textit{RAT} property---the application of \textit{RAT} to the fine-mesh negative discrete Laplacian $A_\f$ does not yield the coarse-mesh operator $A_\c$ obtained from pure geometric multigrid. Using the properties derived above,
\begin{align*}
A_\c &= -D_\c G_\c + \tau_0 E_{0,\c} + \tau_D E_{D,\c} \\
&= -\mathcal{C}(D_\f) \mathcal{C}(G_\f)  + \tau_0 \mathcal{C}(E_{0,\f}) + \tau_D \mathcal{C}(E_{D,\f})
\end{align*}
which differs from the direct coarsening of $A_\f$,
\[
\mathcal{C}(A_\f) = -\mathcal{C}(D_\f G_\f) + \tau_0 \mathcal{C}(E_{0,\f}) + \tau_D \mathcal{C}(E_{D,\f}) \neq A_\c,
\]
since in general $\mathcal{C}(D_\f G_\f) \neq \mathcal{C}(D_\f) \mathcal{C}(G_\f)$. Informally, $\mathcal{C}(D_\f G_\f) u_\c$ interpolates a function $u_\c \in V_{2h}(\mathcal{E}_\c)$ onto the fine mesh $\mathcal{E}_\f$, computes the gradient as a function in $V_h^d(\mathcal{E}_\f)$, computes the divergence as a function in $V_h(\mathcal{E}_\f)$, and projects the result back to the coarse mesh $\mathcal{E}_\c$. On the other hand, $\mathcal{C}(D_\f) \mathcal{C}(G_\f) u_\c$ projects the computed fine-mesh gradient onto the coarse mesh and then immediately interpolates the result in order to compute the discrete divergence on the fine mesh, before projecting the final result back to the coarse mesh. That is,
\[
\mathcal{C}(D_f) \mathcal{C}(G_f) = \mathcal{C}(D_f I_c^f R_f^c G_f) \neq \mathcal{C}(D_f G_f),
\]
since $I_c^f R_f^c \neq \mathbb{I}$.

\subsubsection{Flux coarsening}
The coarse operator $A_\c$ obtained from pure geometric multigrid may be viewed as applying \textit{RAT} to the equations $\vec{q} = \nabla u$ and $-\nabla \cdot \vec{q} = f$ separately. The flux formulation of LDG, \eqref{eq:LDG_block}, naturally displays this coarsening strategy. To show this, note that we can write the flux formulation with input and output in the user-defined basis by premultiplying \eqref{eq:LDG_block} by the inverse mass matrix to obtain
\begin{equation}\label{eq:LDG_block_nodal}
\begin{bmatrix}
\phantom{-}I & -G \\[0.3em] -D & \phantom{-}T
\end{bmatrix}
\begin{bmatrix}
\vec{q}_h \\[0.3em] u_h
\end{bmatrix}
=
\begin{bmatrix}
M^{-1} j \\[0.3em] M^{-1} k
\end{bmatrix}.
\end{equation}
Applying \textit{RAT} in a block fashion to the flux formulation \eqref{eq:LDG_block_nodal} then yields the discrete operator 
\begin{equation}\label{eq:LDG_block_RAT}
\begin{bmatrix}
R_\f^\c & 0 \\[0.3em] 0 & R_\f^\c
\end{bmatrix}
\begin{bmatrix}
I & -G_\f \\[0.3em] -D_\f & \phantom{-}T_\f
\end{bmatrix}
\begin{bmatrix}
I_\c^\f & 0 \\[0.3em] 0 & I_\c^\f
\end{bmatrix}
=
\begin{bmatrix}
I & -G_c \\[0.3em] -D_\c & \phantom{-}T_\c
\end{bmatrix}.
\end{equation}
Taking the Schur complement of the right-hand side of \eqref{eq:LDG_block_RAT}, we obtain
\begin{equation}
\begin{aligned}
A_c &= -D_c G_c + T_c \\
&= -\mathcal{C}(D_f) \mathcal{C}(G_f) + \tau_0 \mathcal{C}(E_{0,f}) + \tau_D \mathcal{C}(E_{D,f}),
\end{aligned}
\end{equation}
which is exactly the coarse operator from pure geometric multigrid. Thus, applying operator coarsening to the flux formulation of LDG, which is equivalent to separately coarsening the equations $\vec{q} = \nabla u$ and $-\nabla \cdot \vec{q} = f$, is the same as pure geometric multigrid.

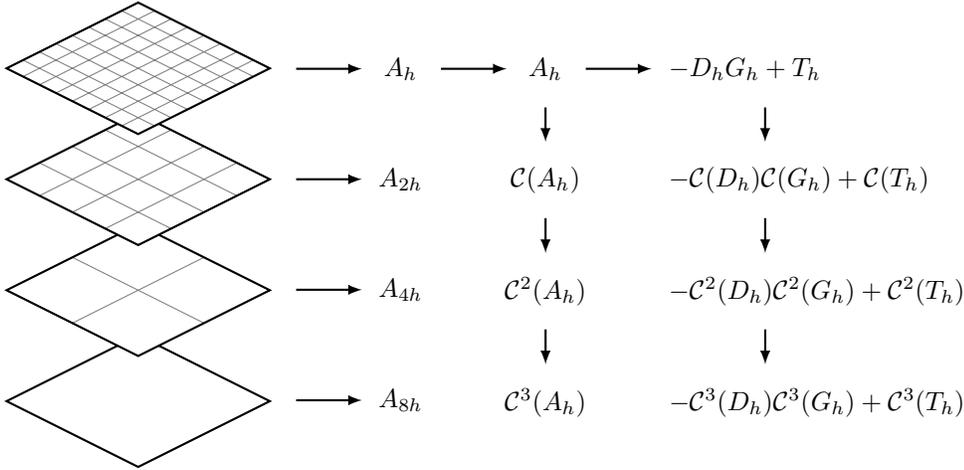
\begin{figure}[t]
\centering
\begin{tikzpicture}[scale=0.35,every node/.style={minimum size=1cm},on grid]

	 \begin{scope}[
    	yshift=-240,every node/.append style={
    	    yslant=0.5,xslant=-1},yslant=0.5,xslant=-1]
        \fill[white,fill opacity=1] (0,0) rectangle (5,5);
        \draw[black, thick] (0,0) rectangle (5,5);
    \end{scope}

	 \begin{scope}[
    	yshift=-120,every node/.append style={
    	    yslant=0.5,xslant=-1},yslant=0.5,xslant=-1]
        \fill[white,fill opacity=1] (0,0) rectangle (5,5);
        \draw[step=25mm, black!50] (0,0) grid (5,5);
        \draw[black, thick] (0,0) rectangle (5,5);
    \end{scope}

	 \begin{scope}[
    	yshift=0,every node/.append style={
    	    yslant=0.5,xslant=-1},yslant=0.5,xslant=-1]
        \fill[white,fill opacity=1] (0,0) rectangle (5,5);
        \draw[step=12.5mm, black!50] (0,0) grid (5,5);
        \draw[black, thick] (0,0) rectangle (5,5);
    \end{scope}
	
    \begin{scope}[
    	yshift=120,every node/.append style={
    	    yslant=0.5,xslant=-1},yslant=0.5,xslant=-1]
        \fill[white,fill opacity=1] (0,0) rectangle (5,5);
        \draw[step=6.25mm, black!50] (0,0) grid (5,5);
        \draw[black, thick] (0,0) rectangle (5,5);
    \end{scope}

    \draw[-latex,thick] (6,6.7) to (8.5,6.7) node[right]{$A_h$};
    \draw[-latex,thick] (6,2.5) to (8.5,2.5) node[right]{$A_{2h}$};
    \draw[-latex,thick] (6,-1.7) to (8.5,-1.7) node[right]{$A_{4h}$};
    \draw[-latex,thick] (6,-5.9) to (8.5,-5.9) node[right]{$A_{8h}$};

    \draw[-latex,thick] (11.5,6.7) to (14,6.7) node[right](Ah) {$A_h$};
    \node[below=1.48 of Ah](CAh) {$\mathcal{C}(A_h)$};
    \draw[-latex,thick] (Ah) to (CAh);
    \node[below=1.48 of CAh](CCAh) {$\mathcal{C}^2(A_h)$};
    \draw[-latex,thick] (CAh) to (CCAh);
    \node[below=1.48 of CCAh](CCCAh) {$\mathcal{C}^3(A_h)$};
    \draw[-latex,thick] (CCAh) to (CCCAh);
    
		\draw[-latex,thick] (17,6.7) to (19.5,6.7) node[right=0.1,text width=3.96cm](DGh) {$-D_h G_h + T_h$};
    \node[below=1.48 of DGh,text width=3.96cm,align=left](CDGh) {$-\mathcal{C}(D_h) \mathcal{C}(G_h) + \mathcal{C}(T_h)$};
    \draw[-latex,thick,transform canvas={xshift=-2em}] (DGh) to (CDGh);
    \node[below=1.48 of CDGh,text width=3.96cm](CCDGh) {$-\mathcal{C}^2(D_h) \mathcal{C}^2(G_h) + \mathcal{C}^2(T_h)$};
    \draw[-latex,thick,transform canvas={xshift=-2em}] (CDGh) to (CCDGh);
    \node[below=1.48 of CCDGh,text width=3.96cm](CCCDGh) {$-\mathcal{C}^3(D_h) \mathcal{C}^3(G_h) + \mathcal{C}^3(T_h)$};
    \draw[-latex,thick,transform canvas={xshift=-2em}] (CCDGh) to (CCCDGh);
\end{tikzpicture}
\caption{Three coarsening methods can be used to generate a hierarchy of operators for multigrid: (left column) pure geometric multigrid, where the coarse operators are defined directly from the corresponding coarse meshes, (center column) primal coarsening, where the coarse operators are defined by applying \textit{RAT} to the fine-mesh discrete Laplacian, and (right column) flux coarsening, where the coarse operators are defined by applying \textit{RAT} to the fine-mesh discrete divergence, discrete gradient, and discrete penalty operators, and recombining the results.\vspace{-2em}}
\label{fig:coarsening_types}
\end{figure}

Figure~\ref{fig:coarsening_types} depicts the three types of coarsening that can be performed, given a hierarchy of meshes. In the left column, pure geometric multigrid defines the coarse operators directly from the coarse meshes; in the center column, primal coarsening applies \textit{RAT} to the discrete Laplacian operator; and in the right column, flux coarsening applies \textit{RAT} separately to the discrete divergence and gradient operators. In the above, we have shown the equivalence of the left and right columns. An implementation of constructing the operator hierarchy using flux coarsening is outlined in Algorithm~\ref{alg:build}.

\subsubsection{Benefits of operator coarsening}
It can be useful to define coarse operators directly from fine operators (e.g., by using \textit{RAT}) rather than via discretizations computed directly from coarse meshes. Since coarse mass matrices can be computed automatically according to~\eqref{eq:coarse_mass}, quadrature schemes do not need to be computed for coarse elements---instead, fine-mesh quadrature rules are coarsened automatically via \eqref{eq:coarse_mass}. Similarly, coarse lifting matrices are not explicitly needed since their contribution to the discrete gradient is automatically computed via $G_c = \mathcal{C}(G_\f)$, and so quadrature rules for coarse faces also do not need to be defined. For implicitly defined meshes, such as the ones shown in section~\ref{sec:implicit_mesh}, computing coarse quadrature rules can be computationally intricate or taxing; the fact that operator coarsening obviates the need for this is a substantial benefit. Additionally, operator coarsening can be efficiently implemented using basic linear algebra operations, e.g., block-sparse matrix multiplication, for which highly optimized and parallelized libraries exist; in contrast, computing discretizations directly from coarse meshes relies heavily on the efficiency of one's own code. It is worth noting that the complexity of constructing the operator hierarchy in Algorithm~\ref{alg:build} is the same as the complexity of the multigrid V-cycle in Algorithm~\ref{alg:vcycle} (i.e., $\mathcal{O}(N)$ for $N$ elements); as an approximate indication, in practice the former takes the same computing time as about three to four applications of a V-cycle.

\subsubsection{Relation to other DG methods}
Although we have focused on the LDG method in our discussion, we expect that other DG methods may require similar care in coarsening fine-grid operators such that they are consistent with a pure geometric multigrid method. For instance, methods for which the numerical flux $\widehat{\vec{q}}_h$ depends on the discrete gradient of $u_h$---such as the BR1~\cite{Bassi_97_01} or Brezzi~\cite{Brezzi_00_01} methods---may need similar treatment, as the contribution from the lifting operator $L$ must be coarsened separately.

Other methods, such as the symmetric interior penalty (SIP) method~\cite{Douglas_76_01, Arnold_82_01}, do not require the discrete divergence and discrete gradient operators to be coarsened separately. To demonstrate this for SIP, start with its corresponding bilinear form for a pure Neumann problem: find $u \in V_h$ such that $a(u,v) = l(v)$ for all $v \in V_h$, where
\[
a(u,v) = (\nabla_h u, \nabla_h v) - \int_{\Gamma_0} (\avg{\nabla_h u} \cdot \jump{v} + \jump{u} \cdot \avg{\nabla_h v}) + \tau \int_{\Gamma_0} \jump{u} \cdot \jump{v}
\]
and
\[
l(v) = (f, v) + \int_{\Gamma_N} h\,v^-,
\]
with $\tau$ scaling inversely to the element size. Now consider a pure geometric multigrid method. Let $u_\c, v_\c \in V_{2h}(\mathcal{E}_\c)$. Then
\begin{align*}
  a_\c(u_\c,v_\c) &= (\nabla_{2h} u_\c, \nabla_{2h} v_\c) - \int_{\Gamma_{0,\c}} (\avg{\nabla_{2h} u_\c} \cdot \jump{v_\c} + \jump{u_\c} \cdot \avg{\nabla_{2h} v_\c}) \\
    &\quad\;+ \tau \int_{\Gamma_{0,\c}} \jump{u_\c} \cdot \jump{v_\c} \\
	&= (\nabla_h I_\c^\f u_\c, \nabla_h I_\c^\f v_\c) - \int_{\Gamma_{0,\f}} (\avg{\nabla_h I_\c^\f u_\c} \cdot \jump{I_\c^\f v_\c} + \jump{I_\c^\f u_\c} \cdot \avg{\nabla_h I_\c^\f v_\c}) \\
	&\quad\;+ \tau \int_{\Gamma_{0,\f}} \jump{I_\c^\f u_\c} \cdot \jump{I_\c^\f v_\c} \\[0.3em]
	&= a_\f(I_\c^\f u_\c, I_\c^\f v_\c),
\end{align*}
where equality holds between the first and second lines because $I_\c^\f$ does not introduce nonzero jumps on new mesh faces. Therefore, as matrices (mapping vectors in the user-defined basis to vectors in the same basis),
\[
u_\c^\trans M_\c A_\c v_\c = (I_\c^\f u_\c)^\trans M_\f A_\f (I_\c^\f v_\c).
\]
Assuming the quadratic form is nondegenerate (which is true since $A_\c$ and $A_\f$ are symmetric positive definite, ignoring the trivial kernel), this implies
\[
A_\c = M_\c^{-1} (I_\c^\f)^\trans M_\f A_\f I_\c^\f = R_\f^\c A_\f I_\c^\f.
\]
This is \textit{RAT} applied to $A_\f$, and so applying pure geometric multigrid to SIP is the same as recursively applying standard (primal) operator coarsening to $A_f$.

\subsection{Multigrid preconditioned conjugate gradient}
Recall that a geometric multigrid method utilizes a combination of relaxation/smoothing together with interpolated approximate solutions of coarsened problems. In the present case, the coarsened problem solves for the correction in a residual equation for the same elliptic problem except on a coarser mesh. In particular, it is important to note that, instead of solving $-\Delta_h u = f$, where $\Delta_h$ is the discrete Laplacian in the chosen basis, one instead solves $-M\Delta_h u = Mf$, as the latter system is symmetric positive (semi)definite. Thus, to appropriately define the coarse mesh problem, one may: (i) calculate the residual of the fine mesh linear system $A_\f x_\f = b_\f$, (ii) multiply the residual by $M_\f^{-1}$ to correctly determine the residual as a piecewise polynomial function, (iii) restrict the residual to the coarse mesh, and then (iv) multiply this residual by $M_\c$ of the coarse mesh. Thus, the coarse mesh problem consists of (approximately) solving for $x_\c$ such that
\[
A_\c x_\c = M_\c (R_\f^\c [M_\f^{-1}(b_\f - A_\f x_\f)]),
\]
which, according to the derived restriction operator \eqref{eq:rmatrix}, conveniently simplifies to
\[
A_\c x_\c = (I_\c^\f)^\trans (b_\f - A_\f x_\f),
\]
and so it is unnecessary to multiply by mass matrices in the implementation of the multigrid method; instead, one can simply apply the transpose of the interpolation matrix. With this consideration in mind, the design of a multigrid V-cycle is relatively straightforward and is outlined in Algorithm~\ref{alg:vcycle}.

\begin{figure}[t]
\begin{center}
\begin{minipage}[t]{0.48\textwidth}
\begin{algorithm}[H]
\caption{\small Construction of coarse operators, $\operatorname{Build}(\mathcal{E}_\f,M_\f,G_\f,T_\f)$}
\begin{algorithmic}[1]
\Require{Fine-mesh operators $M_\f$, $G_\f$, $T_\f$}
\Ensure{List of coarse operators}
\State{$\vec{A} := \{\}$}
\State{$A_\f := G_\f^\trans M_\f G_\f + T_\f$}
\If{$\mathcal{E}_\f$ is not the coarsest mesh}
  \State{$M_\c := (I_\c^\f)^\trans M_\f I_\c^\f$}
  \State{$G_\c := M_\c^{-1} (I_\c^\f)^\trans M_\f G_\f I_\c^\f$}
  \State{$T_\c := (I_\c^\f)^\trans T_\f I_\c^\f$}
  \State{$\vec{A} := \operatorname{Build}(\mathcal{E}_\c,M_\c,G_\c,T_\c)$}
\EndIf
\State{\Return $\{A_\f, \vec{A}\}$}
\end{algorithmic}
\label{alg:build}
\end{algorithm}
\end{minipage}\hfill%
\begin{minipage}[t]{0.48\textwidth}
\begin{algorithm}[H]
\caption{\small Multigrid V-cycle $V(\mathcal{E}_\f,x_\f,b_\f)$ on mesh $\mathcal{E}_\f$ with $\nu$ pre- and post-smoothing steps}
\begin{algorithmic}[1]
\If{$\mathcal{E}_\f$ is the bottom level}
  \State{Solve $A_\f x_\f = b_\f$ directly}
\Else
  \State{Relax $\nu$ times}
  \State{$r_\c := (I_\c^\f)^\trans(b_\f - A_\f x_\f)$}
  \State{$x_\c := V(\mathcal{E}_\c,0, r_\c)$}
  \State{$x_\f := x_\f + I_\c^\f x_\c$}
  \State{Relax $\nu$ times}
\EndIf
\State{\Return $x_\f$}
\end{algorithmic}
\label{alg:vcycle}
\vspace{0.28em}
\end{algorithm}
\end{minipage}
\end{center}\vspace{-1em}
\end{figure}

The V-cycle is designed to preserve the symmetric positive (semi)definite property of the discrete problem, making it suitable for preconditioning the conjugate gradient method. To that end, the relaxation sweeps are performed in a symmetric fashion; for order-dependent relaxation methods such as Gauss--Seidel, the first set of relaxation sweeps uses a given ordering of the unknowns and the second set uses the reverse of that ordering.\footnote{If a relaxation scheme is used that is itself symmetric then there is no need to reverse the ordering of unknowns between pre- and post-smoothing steps, as the V-cycle will automatically preserve symmetry.} A multigrid preconditioned conjugate gradient method~\cite{Tatebe_93_01} (MGPCG) combines the advantages of both solvers: the multigrid preconditioner is effective in the interior of the domain where the elliptic behavior of the matrix dominates, while the conjugate gradient method effectively treats the remaining eigenmodes, which in turn are largely associated with the (weak) imposition of the boundary conditions (and in the case of multi-phase elliptic interface problems, jump conditions on internal interfaces)~\cite{Sussman_99_01,Saye_17_01}. We use a single multigrid V-cycle as a preconditioner in the conjugate gradient method.

\section{Numerical results}\label{sec:results}
In this section, we present numerical experiments to assess the efficacy of flux coarsening for LDG discretizations of elliptic PDEs. As the smoother/relaxation method, we use a block Gauss--Seidel smoother with $\nu = 3$ pre- and post-smoothing steps\footnote{As is typical in multigrid methods, increasing the number of pre- and post-smoothing steps can increase the speed of convergence, i.e., decrease $\rho$ as measured by \eqref{eq:rho}; however, doing so comes at the cost of a more expensive V-cycle and therefore may be less efficient. We observed that $\nu = 3$ gave the best computational efficiency in our numerical experiments in terms of reducing the error by a fixed factor.} in the V-cycle. We initially set the interior penalty parameter to $\tau_0 = 0.01/h$; additional analysis of the influence of penalty parameters is given in section~\ref{sec:tau_study}. We measure multigrid performance via the average convergence factor
\begin{equation} \label{eq:rho}
\rho = \exp\left(\frac{1}{N} \log\frac{\|e_N\|_2}{\|e_0\|_2}\right),
\end{equation}
where $N$ is the number of iterations required to reduce the relative error by a factor of $10^{-10}$ and $e_i$ is the error at iteration $i$ of either standalone multigrid (i.e., $i$ many V-cycles) or MGPCG (i.e., the $i^\text{th}$ iteration of CG preconditioned by a single V-cycle). In effect, $\rho$ measures the average slope of $e_i$ on a log-linear graph. Convergence is measured using a right-hand side of $f=0$ with a random nonzero initial guess for $u$. Convergence results are presented in the following as graphs of $\rho$ as a function of element size $h$, polynomial degree $p$, etc.; the same data is presented in tabular form in the supplementary material attached to this paper.

\subsection{Uniform Cartesian grids}\label{sec:cartesian}
We start by solving \eqref{eq:poisson} with homogeneous Neumann boundary conditions on the domain $\Omega = [0,1]^d$ using a uniform Cartesian grid of size $n \times n$ (for $d=2$) or $n \times n \times n$ (for $d=3$) with cell size $h = 1/n$. We build an $h$-multigrid hierarchy based on uniform grid refinement by applying both primal and flux coarsening to the discretized LDG system, and solve using both standalone V-cycles and MGPCG.

\begin{figure}[p]
  \begin{center}
    \small
    \scalebox{0.9}{\input{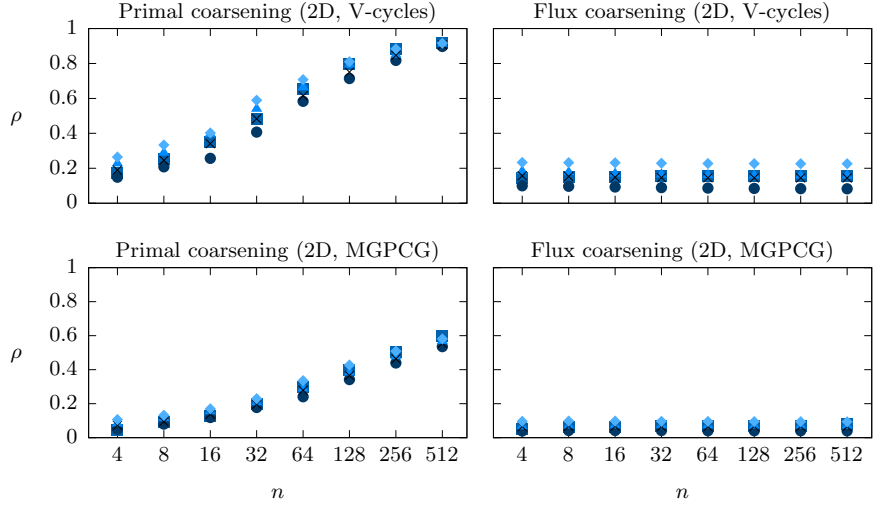}}
  \end{center}
  \caption{$h$-multigrid convergence factors for primal coarsening ($A_\c = \mathcal{C}(A_\f)$) and flux coarsening ($A_\c = -\mathcal{C}(D_\f)\mathcal{C}(G_\f) + \mathcal{C}(T_\f)$) applied to the LDG discretization of Poisson's equation on a uniform $n \times n$ Cartesian grid as ${h\to 0}$. The top row results are computed using V-cycles whereas the bottom row results are computed using MGPCG. The plot markers indicate different polynomial orders: \squaresymbol, \bulletsymbol, \trianglesymbol, \xsymbol, and \diamondsymbol denote $p=1$, $2$, $3$, $4$, and $5$, respectively.}
  \label{fig:uniform_results_2D}
\end{figure}

\begin{figure}[p]
  \begin{center}
    \small
    \scalebox{0.9}{\input{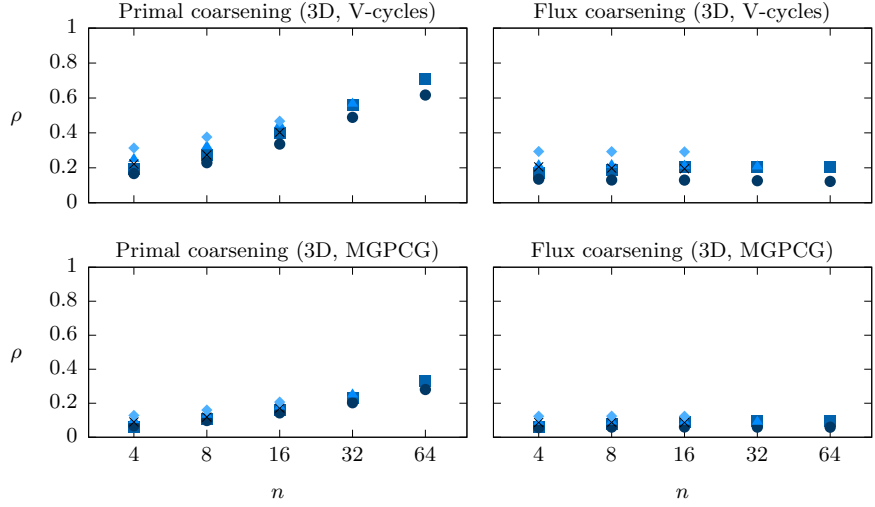}}
  \end{center}
  \caption{$h$-multigrid convergence factors for primal coarsening ($A_\c = \mathcal{C}(A_\f)$) and flux coarsening ($A_\c = -\mathcal{C}(D_\f)\mathcal{C}(G_\f) + \mathcal{C}(T_\f)$) applied to the LDG discretization of Poisson's equation on a uniform $n \times n \times n$ Cartesian grid as ${h\to 0}$. The top row results are computed using V-cycles whereas bottom row results are computed using MGPCG. The plot markers indicate different polynomial orders: \squaresymbol, \bulletsymbol, \trianglesymbol, \xsymbol, and \diamondsymbol denote $p=1$, $2$, $3$, $4$, and $5$, respectively. (Omitted data points correspond to simulations whose memory requirements approximately exceed \SI{120}{\giga\byte}.)}
  \label{fig:uniform_results_3D}
\end{figure}

Figures~\ref{fig:uniform_results_2D} and \ref{fig:uniform_results_3D} show the average convergence factor versus $n$ for polynomial orders $1 \leq p \leq 5$ in 2D and 3D, respectively. In both cases, the multigrid scheme built using flux coarsening exhibits nearly $h$-independent convergence factors of $\rho \approx 0.1$, whereas the scheme based on primal coarsening exhibits poor performance that degrades as $h \to 0$.

\begin{figure}[!t]
  \begin{center}
    \small
    \scalebox{0.9}{\input{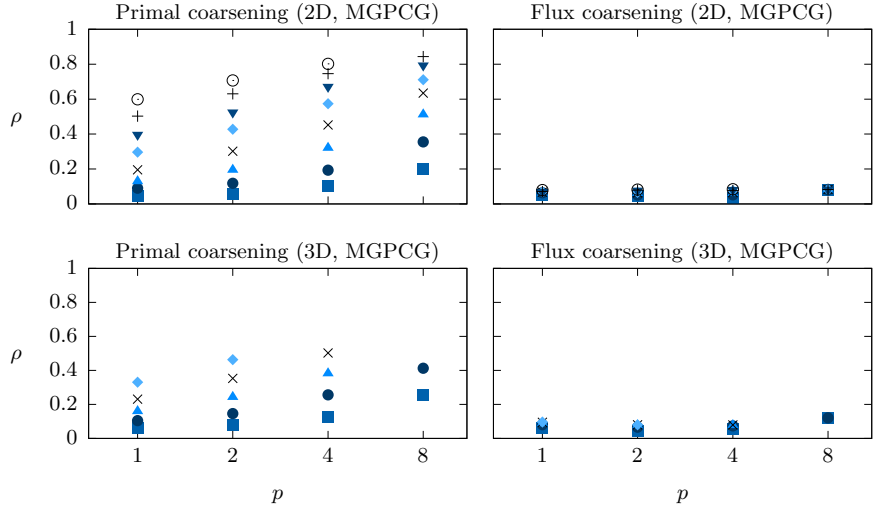}}
  \end{center}
  \caption{$p$-multigrid convergence factors for primal coarsening ($A_\c = \mathcal{C}(A_\f)$) and flux coarsening ($A_\c = -\mathcal{C}(D_\f)\mathcal{C}(G_\f) + \mathcal{C}(T_\f)$) applied to the LDG discretization of Poisson's equation on uniform grids. The $p$-multigrid hierarchy is generated by successively halving the polynomial order (i.e., $p \to p/2 \to p/4 \to \cdots \to 1$). The plot markers indicate different grid sizes: \squaresymbol, \bulletsymbol, \trianglesymbol, \xsymbol, \diamondsymbol, \triangledownsymbol, \plussymbol, and \circlesymbol denote $n=4$, $8$, $16$, $32$, $64$, $128$, $256$, and $512$, respectively. (Omitted data points correspond to simulations whose memory requirements approximately exceed \SI{120}{\giga\byte}.)\vspace{-1em}}
  \label{fig:pmg_results}
\end{figure}

Similar results hold for $p$-multigrid on uniform Cartesian grids. We generate a $p$-multigrid hierarchy by successively halving the polynomial order (i.e., $p \to p/2 \to p/4 \to \cdots \to 1$) and applying both primal and flux coarsening to the discretized LDG system. Figure~\ref{fig:pmg_results} shows convergence factor versus $p$ for grid sizes $4 \leq n \leq 512$ in 2D and 3D. Again, convergence appears to be independent of $p$ (at least up to $p=8$) when flux coarsening is used with $p$-multigrid, whereas performance degrades with increasing $p$ for primal coarsening.

\subsection{On the effect of penalty parameters on multigrid performance}\label{sec:tau_study}

\begin{figure}[!t]
  \centering
  \begin{minipage}{0.49\textwidth}
  \footnotesize
  \scalebox{0.85}{%
  \begin{overpic}[scale=0.35]{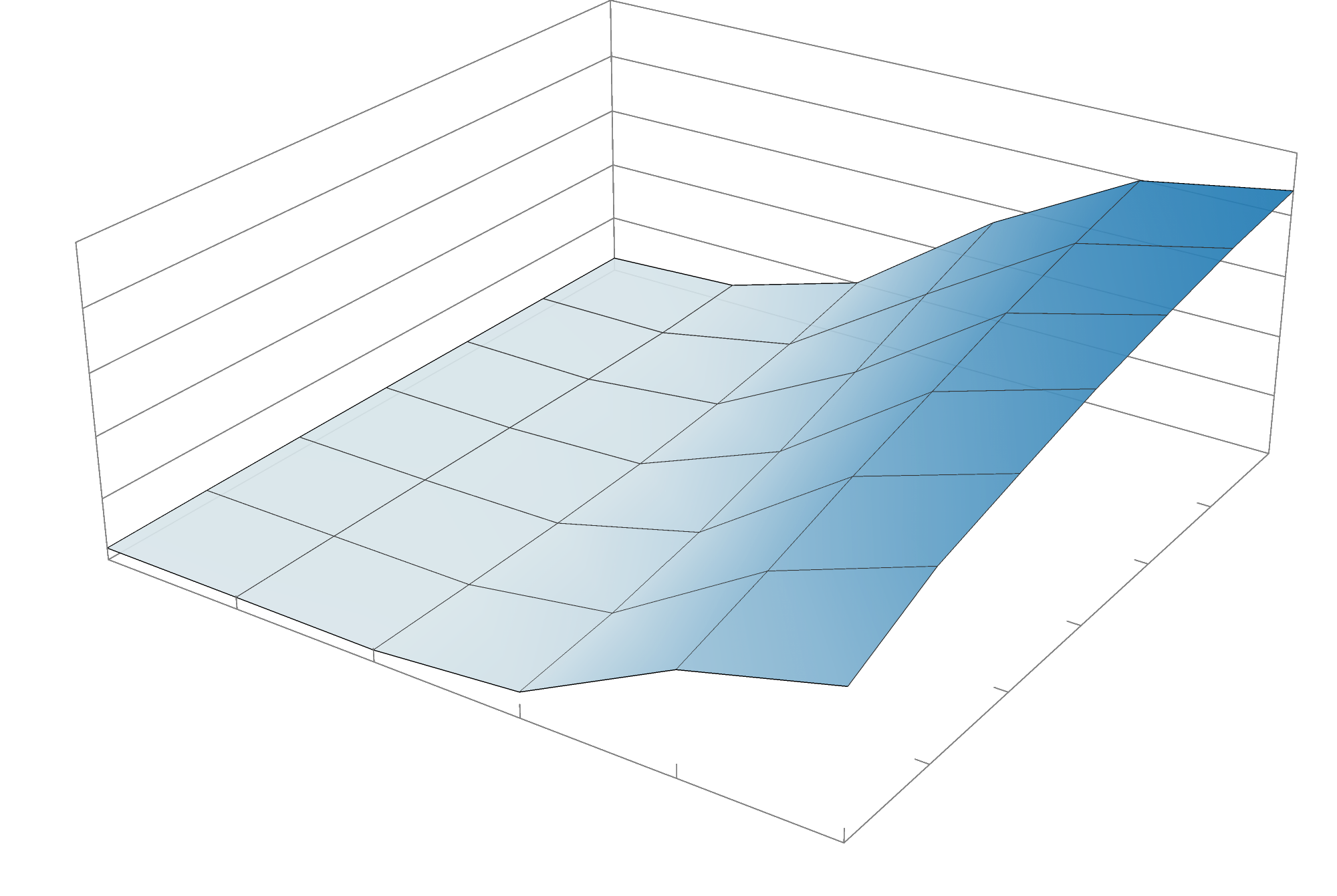}
   \put(-6,35) {$\rho$}
   \put(23,4) {$\widetilde{\tau}_0$}
   \put(85,10) {$n$}
   \put(1.7,47) {$1$}
   \put(-1.2,42) {$0.8$}
   \put(-0.5,37) {$0.6$}
   \put(0,32) {$0.4$}
   \put(0.5,27) {$0.2$}
   \put(4.8,18.8) {$10^{-2}$}
   \put(13.9,15.4) {$10^{-1}$}
   \put(25,11) {$10^0$}
   \put(35,7.5) {$10^1$}
   \put(46,3.4) {$10^2$}
   \put(56,0) {$10^3$}
   \put(65,0.5) {$2^2$}
   \put(70.5,5.5) {$2^3$}
   \put(76,10.5) {$2^4$}
   \put(81.5,15.5) {$2^5$}
   \put(86.5,20.2) {$2^6$}
   \put(91.6,25) {$2^7$}
   \put(96.5,29.8) {$2^8$}
  \end{overpic}%
  }\end{minipage}\hfill%
  \begin{minipage}{0.49\textwidth}
  \footnotesize
  \scalebox{0.85}{%
  \begin{overpic}[scale=0.35]{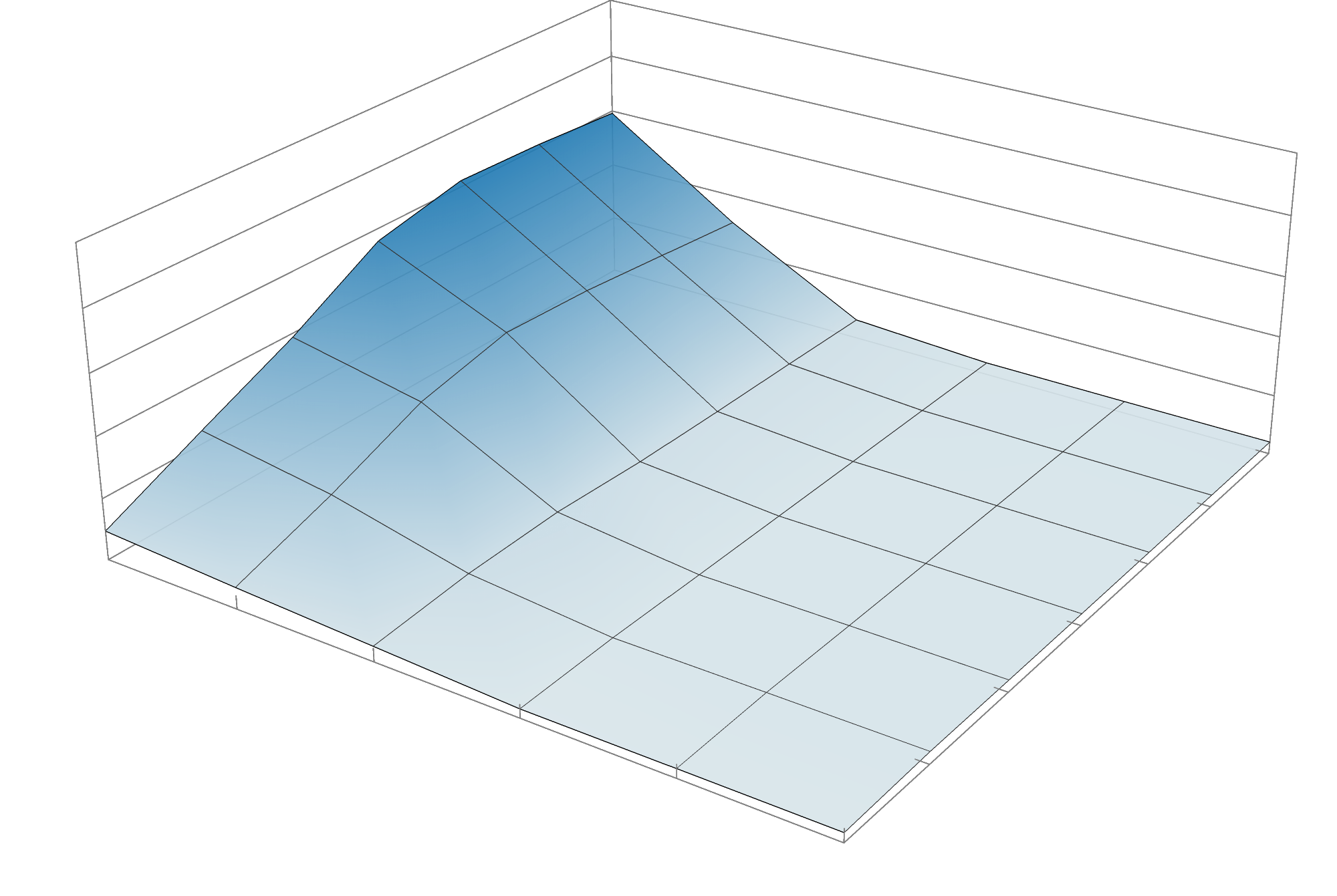}
   \put(-6,35) {$\rho$}
   \put(23,4) {$\widetilde{\tau}_D$}
   \put(85,10) {$n$}
   \put(1.7,47) {$1$}
   \put(-1.2,42) {$0.8$}
   \put(-0.5,37) {$0.6$}
   \put(0,32) {$0.4$}
   \put(0.5,27) {$0.2$}
   \put(4.8,18.8) {$10^{-2}$}
   \put(13.9,15.4) {$10^{-1}$}
   \put(25,11) {$10^0$}
   \put(35,7.5) {$10^1$}
   \put(46,3.4) {$10^2$}
   \put(56,0) {$10^3$}
   \put(65,0.5) {$2^2$}
   \put(70.5,5.5) {$2^3$}
   \put(76,10.5) {$2^4$}
   \put(81.5,15.5) {$2^5$}
   \put(86.5,20.2) {$2^6$}
   \put(91.6,25) {$2^7$}
   \put(96.5,29.8) {$2^8$}
  \end{overpic}%
  }\end{minipage}%
  \caption{Effects of interior penalty parameter $\tau_0 = \widetilde{\tau}_0/h$ (left) and Dirichlet penalty parameter $\tau_D = \widetilde{\tau}_D/h$ (right) on multigrid convergence factor $\rho$ for an $n \times n$ Cartesian mesh with $p=2$, where $h = 1/n$; see section \ref{sec:tau_study}.\vspace{-1em}}
  \label{fig:tau_study}
\end{figure}
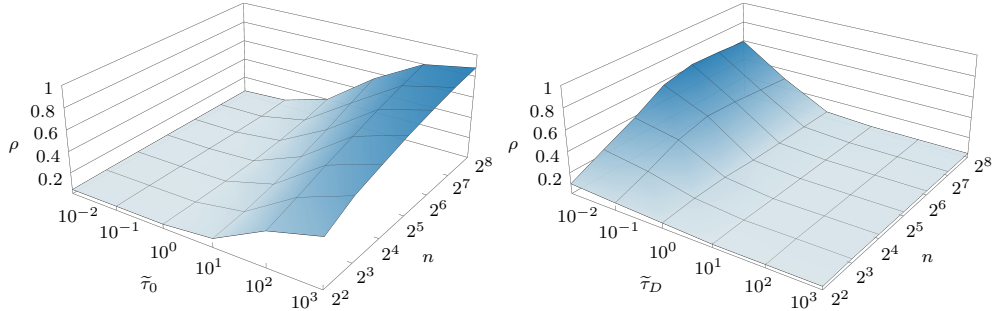

Figure~\ref{fig:tau_study} (left) shows a study of the impact the interior penalty parameter $\tau_0 = \widetilde{\tau}_0/h$ has on multigrid convergence for a Poisson problem on a uniform $n \times n$ mesh with periodic boundary conditions and $p = 2$. Smaller values of $\widetilde{\tau}_0$ yield better convergence factors; for $\widetilde{\tau}_0 > 10$, multigrid performance begins to degrade as the mesh is refined. For the remainder of our tests, we set $\widetilde{\tau}_0 = 0.01$ so that $\tau_0 = 0.01/h$.

In our tests, the imposition or combination of Dirichlet, Neumann, or periodic boundary conditions does not affect the conclusions made in this work. However, for problems with Dirichlet boundary conditions the choice of Dirichlet penalty parameter $\tau_D$ can impact multigrid efficiency. Informally, Dirichlet boundary conditions are enforced in a DG method only weakly and the smoothing/relaxation method of a V-cycle can only effectively enforce the boundary condition if the associated penalty parameter is sufficiently strong. Figure~\ref{fig:tau_study} (right) shows a study of the impact that $\tau_D = \widetilde{\tau}_D/h$ has on multigrid performance for a homogeneous Dirichlet problem on a uniform $n \times n$ mesh. For $\widetilde{\tau}_D < 1$, the average convergence factor $\rho$ degrades as $n$ increases; a good choice in this case appears to be $10<\widetilde{\tau}_D <100$.

Ultimately, the proper choice of both $\tau_0$ and $\tau_D$ is application-dependent and concerns not only multigrid performance but also discretization accuracy and effects of penalty stabilization on the conditioning of the linear systems.

\subsection{Adaptive mesh refinement}\label{sec:amr}
\begin{figure}[!t]
  \centering
  \includegraphics[width=0.8\textwidth]{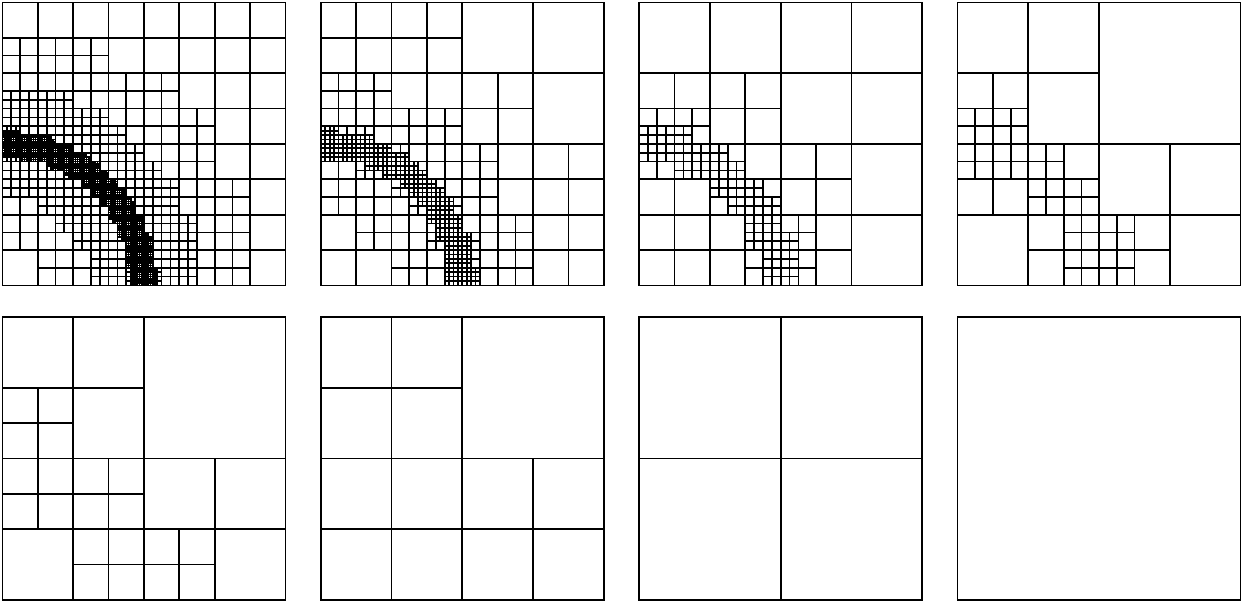}\vspace{-0.5em}
  \caption{An example of a geometric multigrid hierarchy inherited from an adaptively refined quadtree with rapid coarsening. The finest level depicted has smallest cell size equal to $h_\text{min}=1/128$}.
  \label{fig:amr}
\end{figure}

Next, we solve the Neumann problem \eqref{eq:poisson} on an adaptively refined Cartesian mesh, where refinement is performed according to some prescribed spatially-varying function. We implement adaptivity using a quadtree (in 2D) or octree (in 3D), which naturally defines a geometric hierarchy of meshes. An example hierarchy is shown in Figure~\ref{fig:amr}. Note that some elements at various levels of the hierarchy have four neighbors on a single side, so that large elements and small elements may share part of a face.

\begin{figure}[!t]
  \begin{center}
    \small
    \scalebox{0.9}{\input{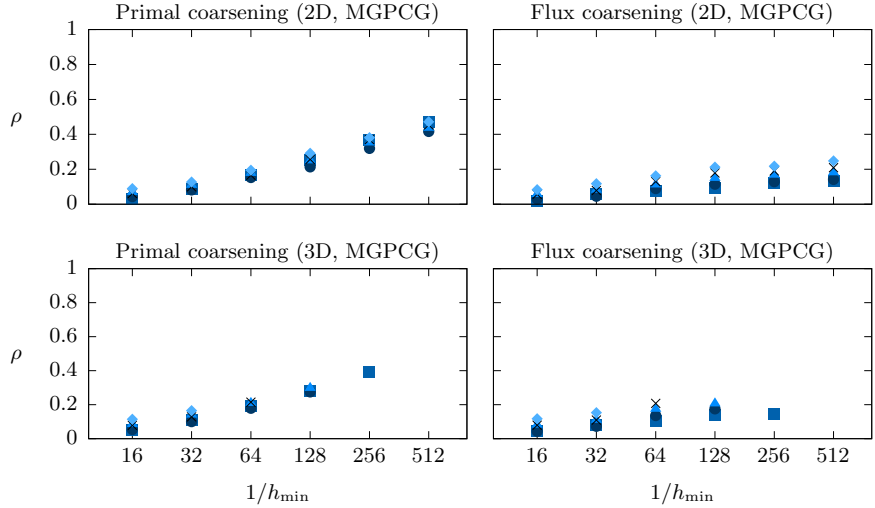}}
  \end{center}
  \caption{$h$-multigrid convergence factors for primal ($A_\c = \mathcal{C}(A_\f)$) and flux ($A_\c = -\mathcal{C}(D_\f)\mathcal{C}(G_\f) + \mathcal{C}(T_\f)$) coarsening applied to the LDG discretization of Poisson's equation on an adaptively refined grid in 2D and 3D. The plot markers indicate different polynomial orders: \squaresymbol, \bulletsymbol, \trianglesymbol, \xsymbol, and \diamondsymbol denote $p=1$, $2$, $3$, $4$, and $5$, respectively. (Omitted data points correspond to simulations whose memory requirements approximately exceed \SI{120}{\giga\byte}.)\vspace{-2em}}
  \label{fig:adaptive_results}
\end{figure}

Figure~\ref{fig:adaptive_results} shows the average convergence factor versus $1/h_\text{min}$ for polynomial orders $1 \leq p \leq 5$, where $h_\text{min}$ is the size of the smallest element in the mesh. In both 2D and 3D, the multigrid method based on primal coarsening exhibits performance that degrades as $h_\text{min} \to 0$, whereas the method based on flux coarsening yields good performance that is nearly independent of $h_\text{min}$ for all $p$ considered.

\subsection{Implicitly defined meshes and elliptic interface problems}\label{sec:implicit_mesh}

Our last two examples are designed to exemplify the benefits of operator coarsening by considering cases in which a pure geometric multigrid method would be intricate or difficult to implement, such as on nontrivial domains with complex geometry or for elliptic interface problems in which the interface has extreme geometry. The first example consists of a curved domain containing holes and thin pieces, and the second example is a multi-phase elliptic interface problem with small circles, filaments, and cusps in the interface geometry. In both cases, we make use of a recently developed framework for computing high-order accurate multi-phase multi-physics using implicitly defined meshes~\cite{Saye_17_01, Saye_17_02}. The framework shares some aspects with cut-cell techniques wherein a level set function defining the domain geometry or internal interfaces is used to cut through the cells of a background quadtree or octree; tiny cut cells are then merged with neighboring cells to create a mesh in which the shapes of interfacial elements are defined implicitly by the level set function. Quadrature rules for curved elements and nontrivial mesh faces are then computed using high-order accurate schemes for computing integrals on implicitly defined domains restricted to hyperrectangles~\cite{Saye_15_01}; these quadrature schemes are then used in the DG weak formulation, e.g., for computing mass matrices and the lifting operator $L$ on the finest-level mesh.

In both examples, we consider an elliptic PDE problem with Dirichlet boundary conditions. The Dirichlet penalty parameter is chosen to scale inversely with $h$, the typical element size on the finest mesh, such that $\tau_D = 100/h$; the value 100 was determined empirically as being approximately the smallest possible while giving good multigrid performance. 
To measure the convergence rate, we apply the MGPCG method to a homogeneous problem with random nonzero initial condition and measure the average convergence rate using \eqref{eq:rho}.

\begin{figure}[!t]
	\centering
	\footnotesize
	\sffamily
	\def\xx{0.75in}
	\begin{tabular}{@{}c@{\hspace{1mm}}c@{\hspace{1mm}}c@{\hspace{1mm}}c@{\hspace{1mm}}c@{\hspace{1mm}}c@{}}
		\includegraphics[angle=90,width=\xx]{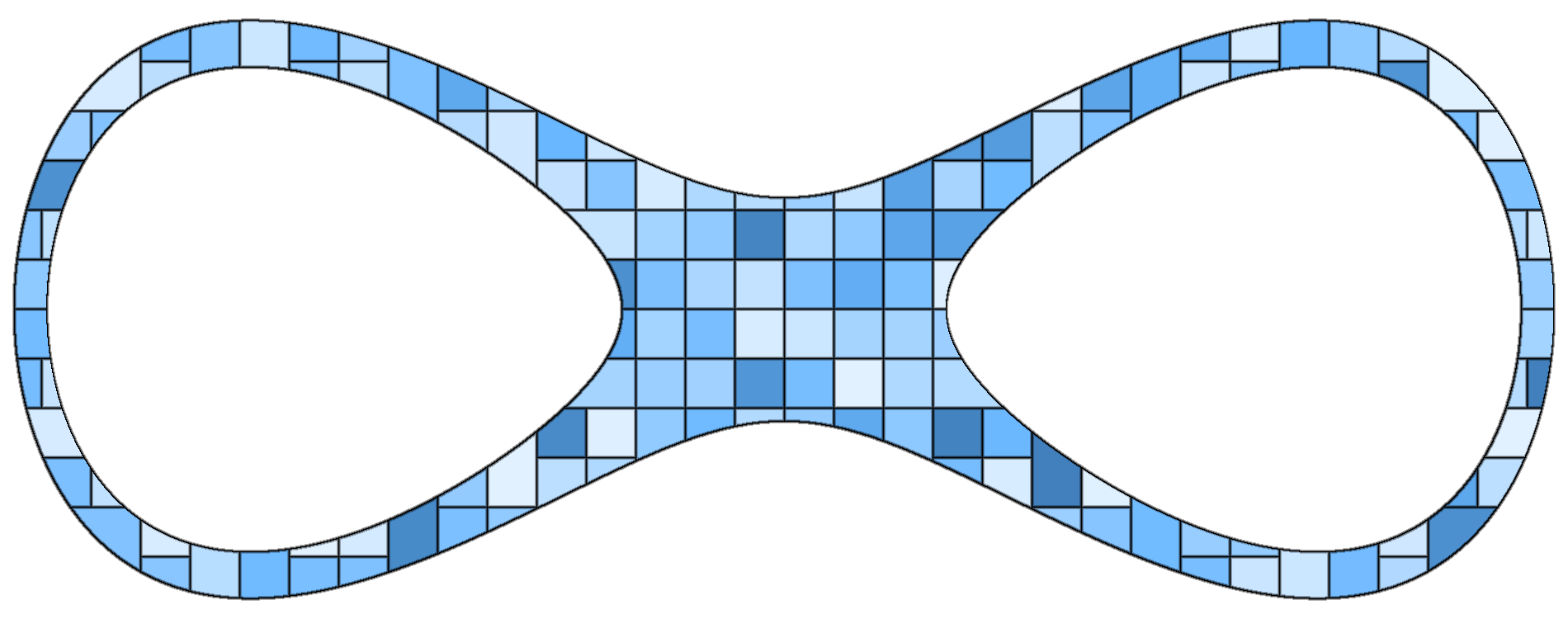} &
		\includegraphics[angle=90,width=\xx]{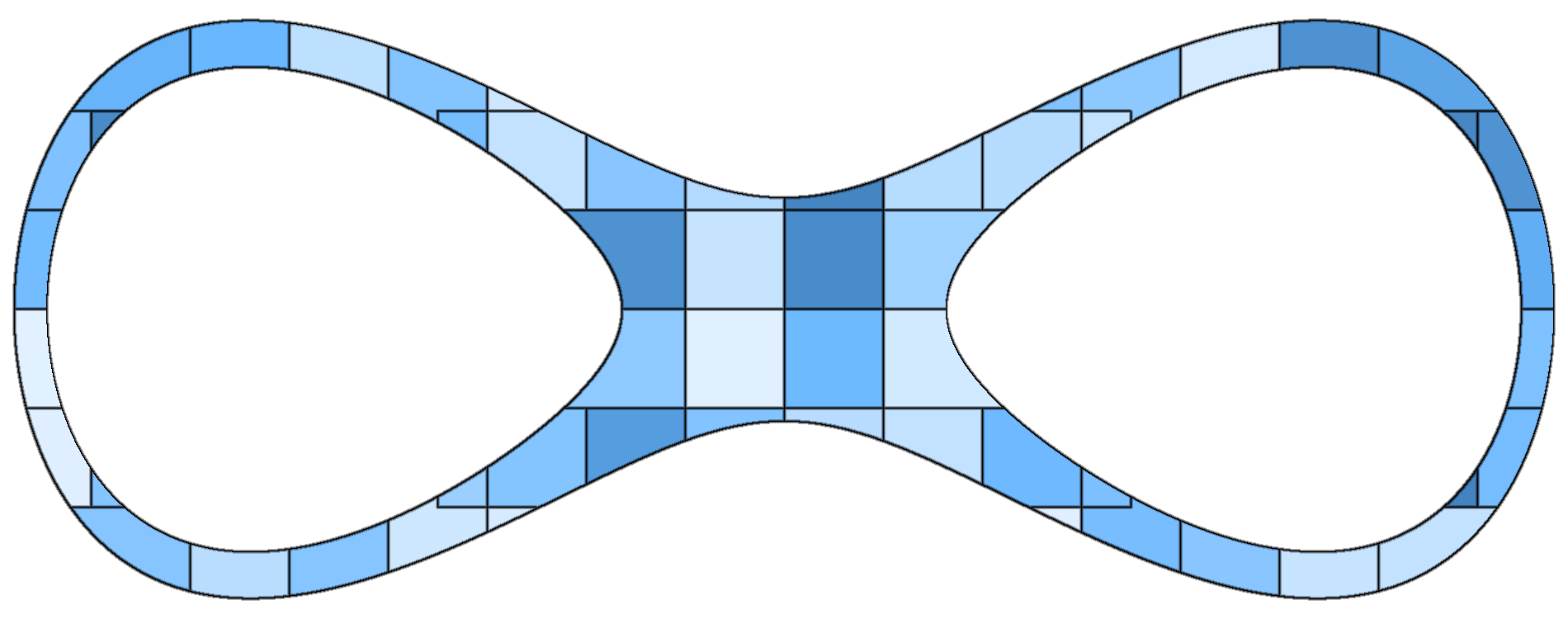} &
		\includegraphics[angle=90,width=\xx]{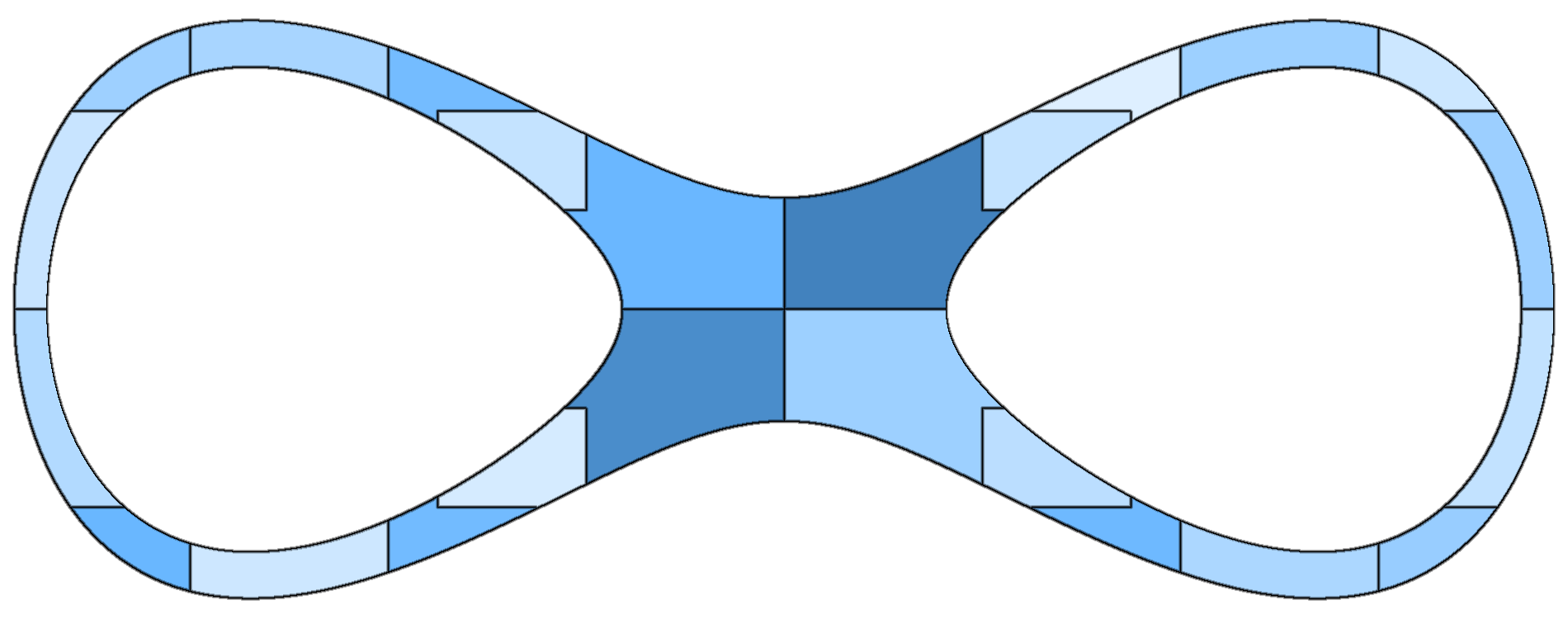} &
		\includegraphics[angle=90,width=\xx]{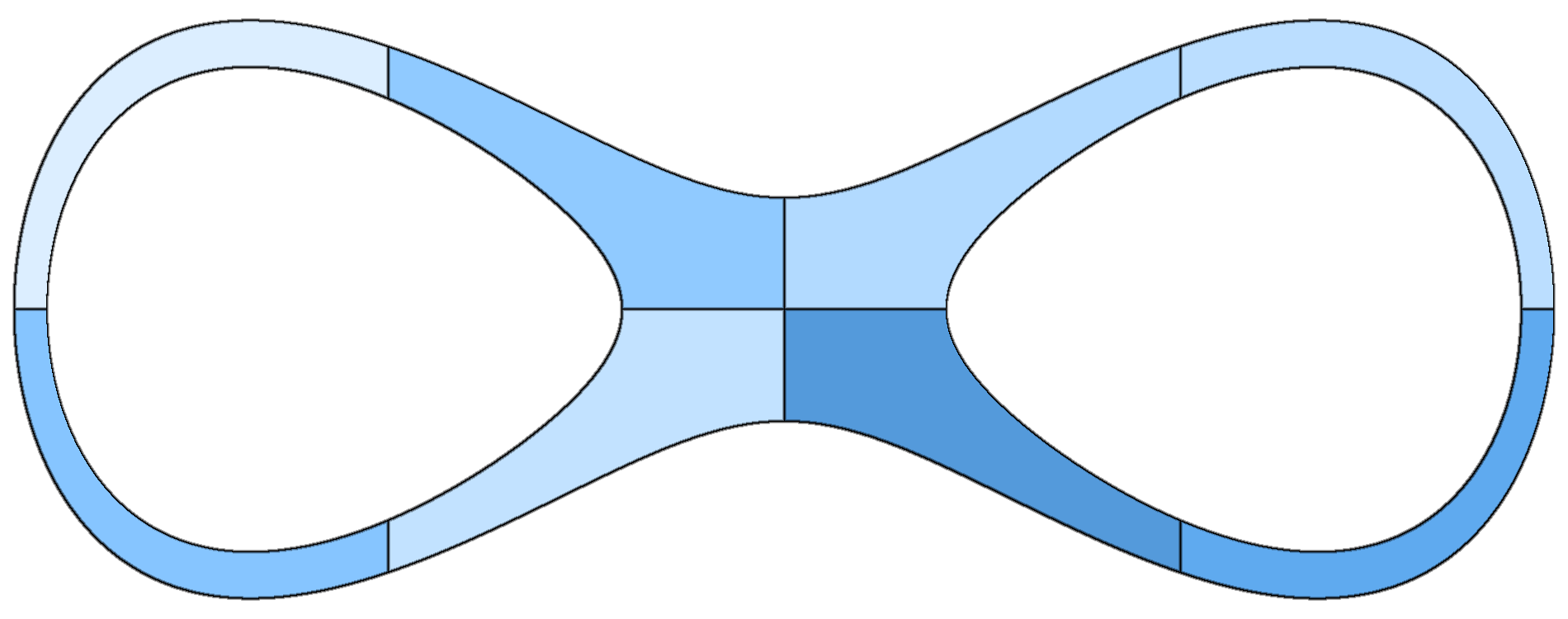} &
		\includegraphics[angle=90,width=\xx]{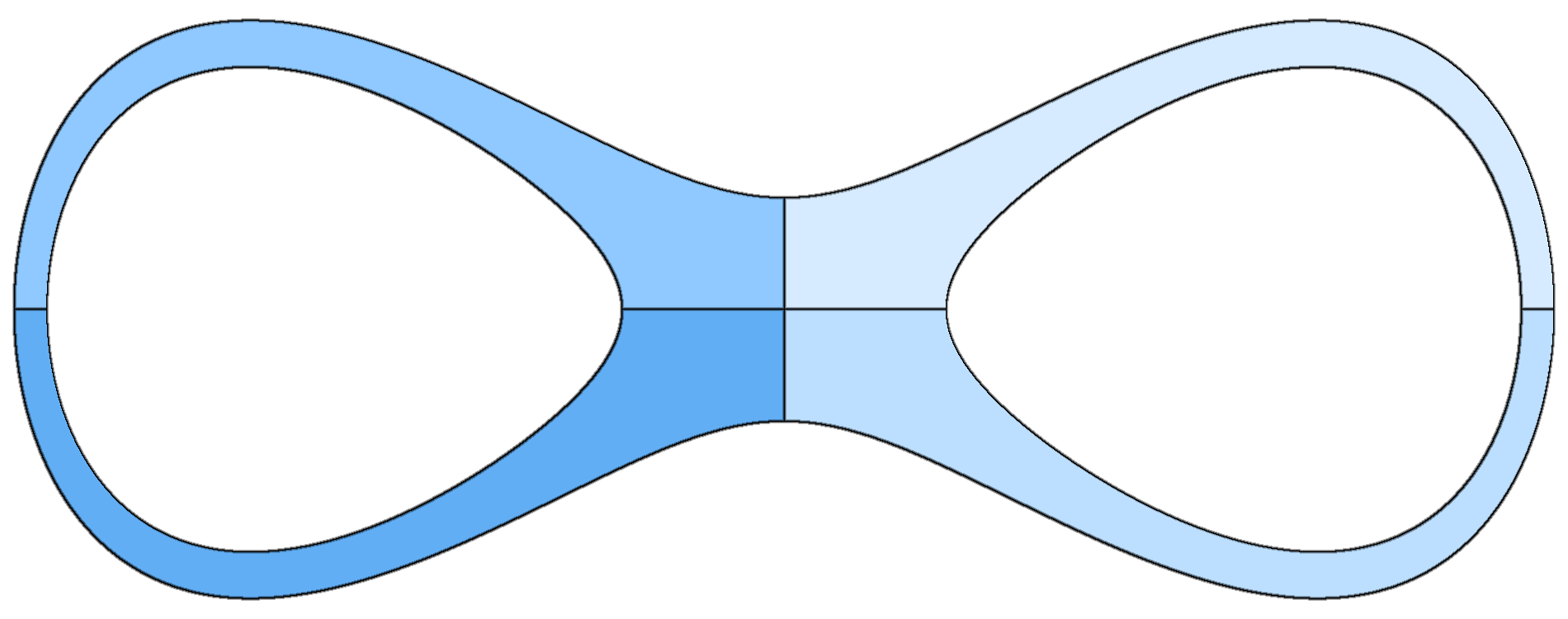} &
		\includegraphics[angle=90,width=\xx]{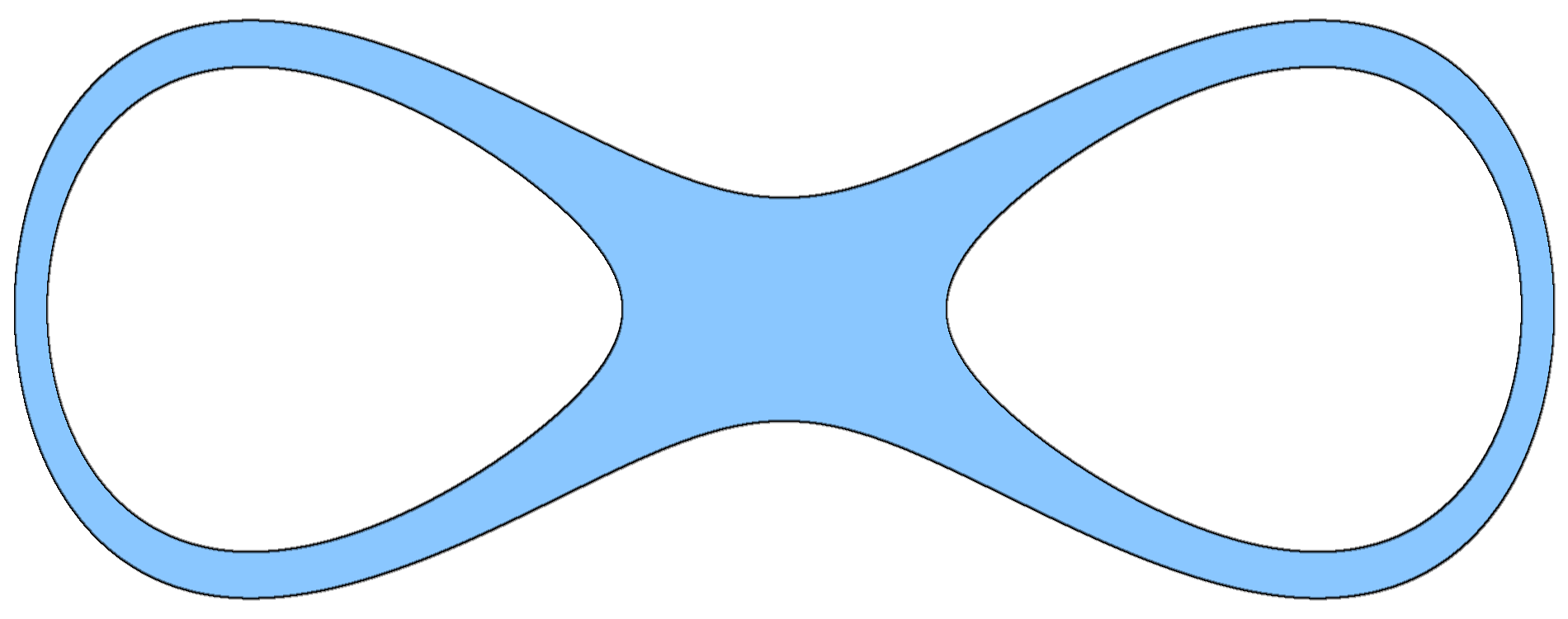} \\
		Fine mesh & & & & & Bottom level
	\end{tabular}\vspace{-0.5em}
	\caption{A single-phase test case applied to a curved domain using implicitly defined meshes. Depicted is the implied $h$-multigrid hierarchy (elements are randomly colored). The coarse meshes are not explicitly constructed in our operator coarsening approach; instead, the discrete gradient and penalty operators and mass matrices are constructed top-down at each level (see Algorithm~\ref{alg:build}). In particular, note that the bottom level of the hierarchy consists of a mesh with a single element containing two holes.\vspace{-2em}}
  \label{fig:lemniscate}
\end{figure}

The first example of a single-phase Poisson problem on a curved domain is illustrated in Figure~\ref{fig:lemniscate} and consists of a figure-eight domain with two holes surrounded by thin segments. It is important to note that the illustrated mesh hierarchy is implicitly formed by our operator coarsening scheme in Algorithm~\ref{alg:build} and it is only the finest-level mesh which is built. On coarser levels, the elements are agglomerated according to the coarsening of the background quadtree; in particular, we note that the bottom level consists of a single element containing two holes. Using the operator coarsening strategy, there is no need to compute quadrature rules for the coarse levels of the mesh hierarchy---the quadrature rules from the fine mesh are effectively coarsened automatically. The results for solving the Dirichlet problem~\eqref{eq:poisson} using flux coarsening\footnote{Results using primal coarsening are similar to previous examples that use primal coarsening---poor multigrid performance is observed that degrades with mesh size---and have been omitted for brevity.} and $h$-multigrid on the curved domain of Figure~\ref{fig:lemniscate} are shown in Figure~\ref{fig:implicit_results} (left); we observe good multigrid convergence factors of $\rho \approx 0.05$--$0.2$, nearly independent of grid size, for $1 \leq p \leq 5$.

\begin{figure}[!t]
	\centering
	\footnotesize
	\sffamily
	\def\svgwidth{0.85\textwidth}
	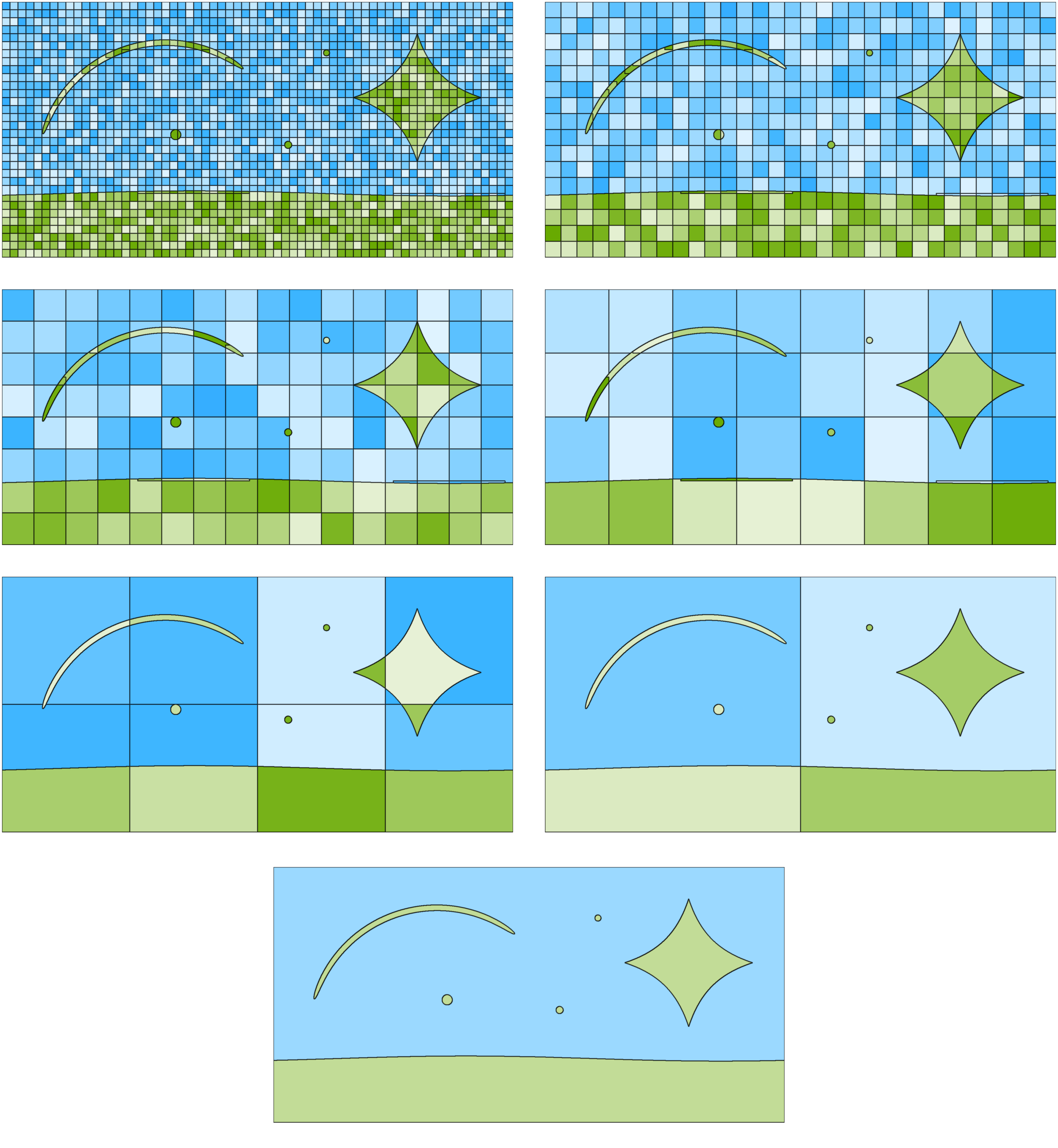%
	\vspace{-0.5em}
	\caption{A test case involving a multi-phase elliptic interface problem. Depicted is the implied $h$-multigrid hierarchy, wherein the elements are randomly colored a shade of green or blue for phase one or two, respectively. Similar to Figure~\ref{fig:lemniscate}, only the finest mesh is explicitly constructed; on coarser levels, the discrete gradient and penalty operators and mass matrices are constructed top-down at each level (see Algorithm~\ref{alg:build}). Note that the bottom level of the hierarchy consists of a mesh of just two elements---the blue element has one component with five holes, whereas the green element consists of six connected components. See the discussion following \eqref{eq:interface} for a description of the noted features \textsf{A}, \textsf{B}, \textsf{C}, \textsf{D\textsubscript{1}}, and \textsf{D\textsubscript{2}}.\vspace{-2em}}
  \label{fig:interface}
\end{figure}

The second example considers a two-phase elliptic interface problem in a rectangular domain illustrated in Figure~\ref{fig:interface}. The corresponding PDE consists of solving
\begin{equation}\label{eq:interface}
\begin{aligned}
-\nabla \cdot (\mu_i \nabla u) &= f &&\text{in } \Omega_i, &&& \jump{u} &= g_\Gamma &&\text{on } \Gamma, \\
u &= g_{\partial} &&\text{on } \partial \Omega, &&& \jump{\mu \nabla u \cdot \vec{n}} &= h_\Gamma &&\text{on } \Gamma,
\end{aligned}
\end{equation}
where $\Gamma$ is the interface between phases $\Omega_1$ (green region, with ellipticity coefficient\footnote{The chosen multi-phase elliptic interface problem has a somewhat mild coefficient jump of a factor of four across the interface. For much larger ratios, e.g., $10^3$ to $10^6$ and beyond, the performance degrades. In these cases, modifications to the LDG discretization can improve accuracy, conditioning, and multigrid performance \cite{Saye_17_01} and will be reported on in forthcoming work.} $\mu_1 = 1$) and $\Omega_2$ (blue region, with $\mu_2 = 4$), and $\jump{\cdot}$ denotes the interfacial jump in the indicated quantity. In this test case, the geometry of the interface has been designed to be challenging---the crescent shape is long and thin; there are three isolated, small circles; and the star shape has sharp cusp-like corners. Once more we note that the illustrated hierarchy in Figure~\ref{fig:interface} is implicitly formed by the operator coarsening strategy, and only the finest-level mesh is actually built. However, we note that the agglomeration strategy for this multi-phase problem is slightly different from the previous test examples---here, elements are only agglomerated with elements belonging to the same phase. Thus, the interface remains sharp throughout all levels of the hierarchy, and this can dramatically improve the performance of multigrid methods for elliptic interface problems, especially when using high-order accurate techniques~\cite{Saye_17_01}. Owing to this agglomeration strategy, coarse mesh levels can have intricate element shapes. Some example features are noted in Figure~\ref{fig:interface}---\textsf{A} indicates a tiny green-phase element surrounded by a large blue-phase element; \textsf{B} indicates two sliver elements; \textsf{C} indicates an element whose cusp-like corners would make it rather difficult to apply a black-box quadrature scheme if one were to explicitly build the coarse-level mesh; and \textsf{D\textsubscript{1,2}} show two green-phase elements, each with multiple connected components (three for \textsf{D\textsubscript{1}} and four for \textsf{D\textsubscript{2}}), which is perhaps rather unusual for a finite element method. Another aspect which motivated the present work on operator coarsening is that it would be nearly impossible to directly apply the cell merging algorithms underlying implicitly defined meshes~\cite{Saye_17_01} to these coarse levels. Although elements with extreme shapes like these (especially tiny elements next to large elements) can traditionally be of concern for numerical discretization of PDEs, according to our tests they pose no problem when present in the coarse levels of a multigrid solver. Results for solving the homogeneous version of \eqref{eq:interface} ($f$, $g_\partial$, $g_\Gamma$, and $h_\Gamma$ all zero) with a random nonzero initial guess using flux operator coarsening are shown in Figure~\ref{fig:implicit_results} (right). In this case of an implicitly defined mesh for which different cell merging decisions take place depending on the refinement of the background grid, leading to different mesh topologies as $n$ is increased, we naturally expect some amount of noise in $\rho$. For the majority of grid sizes, we see that the convergence factor $\rho$ is relatively constant, taking values $\rho \approx 0.15$--$0.3$ reflective of the challenging interface geometry; meanwhile, for the largest mesh corresponding to $n = 1024$, the slight increase in the convergence rate $\rho$ for all $p$ is attributed to the increased ill-conditioning of the system.

\begin{figure}
  \begin{center}
    \small
    \scalebox{0.9}{\input{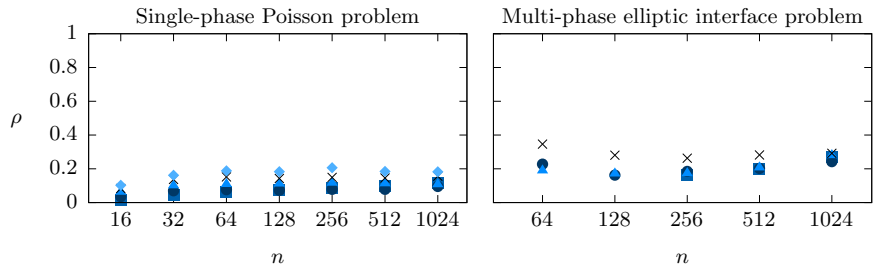}}
  \end{center}\vspace{-1em}
  \caption{$h$-multigrid convergence factors for flux coarsening and MGPCG applied to the LDG discretization of (left) the single-phase Poisson problem on the curved domain shown in Figure~\ref{fig:lemniscate} and (right) the multi-phase elliptic interface test problem in \eqref{eq:interface} on the domain shown in Figure~\ref{fig:interface}. The grid size $n$ is the number of cells of the background Cartesian grid required to cover the longest extent of the domain. The plot markers indicate different polynomial orders: \squaresymbol, \bulletsymbol, \trianglesymbol, \xsymbol, and \diamondsymbol denote $p=1$, $2$, $3$, $4$, and $5$, respectively.\vspace{-2em}}
   \label{fig:implicit_results}
\end{figure}

\section{Concluding remarks}

We have presented an $hp$-multigrid method for LDG discretizations of elliptic problems that is based on coarsening the discrete gradient and divergence operators from the flux formulation. We have shown that coarsening fine-grid operators in this way results in a method that is equivalent to pure geometric multigrid, but avoids the need to compute quantities associated with coarse meshes, such as lifting operators and quadrature rules. Whereas traditional Galerkin operator coarsening applied to the primal formulation exhibits poor multigrid performance, operator coarsening applied to the flux formulation performs well---convergence factors are nearly independent of both mesh size $h$ and polynomial order $p$ for the demonstrated test problems on uniform Cartesian grids, adaptively refined meshes, and implicitly defined meshes on complex geometries. %

Though most of our analysis has focused on the LDG method, we believe that the essential observation applies to other forms of DG discretization of elliptic problems, particularly those in which lifting operators enter the numerical flux for $\vec q$. A more thorough analysis for other DG methods and more general choices of numerical fluxes would be required to determine whether the multigrid method described here extends to other methods, such as CDG or HDG. Similarly, though we have employed equal-order elements in this work, i.e., polynomials of the same degree for both $u$ and $\vec q$, operator-coarsening for mixed-order elements~\cite{Brezzi_05_01} would be an interesting topic for future investigation.

We considered in this work structured meshes (Cartesian, quadtree, and octree meshes) as well as semi-unstructured, nonconforming, implicitly defined meshes that result from cell merging procedures (see Figures \ref{fig:lemniscate} and \ref{fig:interface}). Applying the multigrid ideas presented here to problems involving more general unstructured meshes is currently under investigation. In this setting, it may be worthwhile to consider different types of relaxation methods owing to their critical role in the overall efficacy of a multigrid method. For example, additive Schwarz smoothers have been shown effective on non-nested polygonal meshes resulting from agglomeration procedures~\cite{Antonietti_19_01}; these smoothers could be studied in the flux coarsening context as well. %

The idea of coarsening the divergence and gradient operators separately may also be useful for AMG methods, which currently treat the discrete Laplacian operator in its entirety as a black box. Indeed, black-box AMG algorithms applied to LDG discretizations appear to struggle~\cite{Olson_11_01}; perhaps applying AMG separately to the divergence and gradient operators in the flux formulation may yield better results.

\bibliographystyle{siamplain}
\bibliography{references}

\end{document}


\maketitle

This supplementary material accompanies the article \emph{Efficient Operator-Coarsening Multigrid Schemes for Local Discontinuous Galerkin Methods}, published in the SIAM Journal on Scientific Computing.

\renewcommand{\thesection}{S}
\section{Tables of multigrid convergence factors}
Here we present the raw multigrid convergence factors used to generate Figures 4.1, 4.2, 4.3, 4.6, and 4.9 of the main paper. The data are arranged below in tables, with grid size, polynomial order, coarsening strategy, and solver indicated. Omitted data points correspond to simulations whose memory requirements approximately exceed \SI{120}{\giga\byte}.

\begin{table}[htb]
\caption{$h$-multigrid convergence factors for primal and flux coarsening applied to the LDG discretization of Poisson's equation on a uniform $n \times n$ Cartesian grid as ${h\to 0}$ using V-cycles (see Figure 4.1).\vspace{-1.8em}}
\label{tab:2D_uniform_vcycles}
\centering
\small
\begin{tabular}{@{}cl@{\qquad}cccccccc@{}}
&&&&&&&&& \\ \midrule
&& \multicolumn{8}{c}{$n$} \\
& $p$ & 4 & 8 & 16 & 32 & 64 & 128 & 256 & 512 \\ \cmidrule{1-10}
\multirow{5}{*}{Primal coarsening}
& 1 & 0.17 & 0.25 & 0.35 & 0.48 & 0.65 & 0.80 & 0.88 & 0.92 \\ 
& 2 & 0.15 & 0.21 & 0.26 & 0.41 & 0.58 & 0.71 & 0.82 & 0.90 \\ 
& 3 & 0.23 & 0.29 & 0.40 & 0.54 & 0.66 & 0.79 & 0.85 & 0.91 \\ 
& 4 & 0.19 & 0.25 & 0.34 & 0.48 & 0.62 & 0.75 & 0.85 & 0.91 \\ 
& 5 & 0.26 & 0.33 & 0.40 & 0.59 & 0.71 & 0.81 & 0.88 & 0.92 \\ \cmidrule{1-10}
\multirow{5}{*}{Flux coarsening}
& 1 & 0.14 & 0.15 & 0.15 & 0.16 & 0.15 & 0.15 & 0.15 & 0.15 \\ 
& 2 & 0.10 & 0.10 & 0.09 & 0.09 & 0.09 & 0.08 & 0.08 & 0.08 \\ 
& 3 & 0.18 & 0.18 & 0.17 & 0.16 & 0.16 & 0.16 & 0.16 & 0.16 \\ 
& 4 & 0.16 & 0.15 & 0.15 & 0.15 & 0.15 & 0.14 & 0.14 & 0.14 \\ 
& 5 & 0.23 & 0.23 & 0.23 & 0.23 & 0.23 & 0.23 & 0.23 & 0.23 \\ \midrule
\end{tabular}
\end{table}

\begin{table}[htb]
\caption{$h$-multigrid convergence factors for primal and flux coarsening applied to the LDG discretization of Poisson's equation on a uniform $n \times n$ Cartesian grid as ${h\to 0}$ using MGPCG (see Figure 4.1).\vspace{-1.8em}}
\label{tab:2D_uniform_mgpcg}
\centering
\small
\begin{tabular}{@{}cl@{\qquad}cccccccc@{}}
&&&&&&&&& \\ \midrule
&& \multicolumn{8}{c}{$n$} \\
& $p$ & 4 & 8 & 16 & 32 & 64 & 128 & 256 & 512 \\ \cmidrule{1-10}
\multirow{5}{*}{Primal coarsening}
& 1 & 0.04 & 0.09 & 0.13 & 0.20 & 0.30 & 0.40 & 0.50 & 0.60 \\ 
& 2 & 0.05 & 0.08 & 0.12 & 0.18 & 0.24 & 0.34 & 0.44 & 0.54 \\ 
& 3 & 0.08 & 0.11 & 0.15 & 0.21 & 0.29 & 0.38 & 0.48 & 0.57 \\ 
& 4 & 0.07 & 0.09 & 0.14 & 0.20 & 0.28 & 0.37 & 0.46 & 0.56 \\ 
& 5 & 0.11 & 0.13 & 0.17 & 0.23 & 0.34 & 0.43 & 0.51 & 0.58 \\ \cmidrule{1-10}
\multirow{5}{*}{Flux coarsening}
& 1 & 0.05 & 0.06 & 0.06 & 0.07 & 0.07 & 0.07 & 0.07 & 0.08 \\ 
& 2 & 0.04 & 0.04 & 0.04 & 0.04 & 0.04 & 0.04 & 0.04 & 0.04 \\ 
& 3 & 0.07 & 0.07 & 0.07 & 0.07 & 0.07 & 0.07 & 0.07 & 0.07 \\ 
& 4 & 0.07 & 0.07 & 0.07 & 0.07 & 0.06 & 0.07 & 0.07 & 0.07 \\ 
& 5 & 0.10 & 0.10 & 0.10 & 0.10 & 0.10 & 0.09 & 0.09 & 0.09 \\ \midrule
\end{tabular}
\end{table}

\begin{table}[htb]
\caption{$h$-multigrid convergence factors for primal and flux coarsening applied to the LDG discretization of Poisson's equation on a uniform $n \times n \times n$ Cartesian grid as ${h\to 0}$ using V-cycles (see Figure 4.2).\vspace{-1.8em}}
\label{tab:3D_uniform_vcycles}
\centering
\small
\begin{tabular}{@{}cl@{\qquad}ccccc@{}}
&&&&&& \\ \midrule
&& \multicolumn{5}{c}{$n$} \\
& $p$ & 4 & 8 & 16 & 32 & 64 \\ \cmidrule{1-7}
\multirow{5}{*}{Primal coarsening}
& 1 & 0.19 & 0.27 & 0.40 & 0.56 & 0.71 \\ 
& 2 & 0.17 & 0.23 & 0.34 & 0.49 & 0.62 \\ 
& 3 & 0.25 & 0.32 & 0.45 & 0.57 &  --  \\ 
& 4 & 0.22 & 0.27 & 0.40 &  --  &  --  \\ 
& 5 & 0.31 & 0.38 & 0.47 &  --  &  --  \\ \cmidrule{1-7}
\multirow{5}{*}{Flux coarsening}
& 1 & 0.17 & 0.19 & 0.20 & 0.20 & 0.21 \\ 
& 2 & 0.14 & 0.13 & 0.13 & 0.13 & 0.12 \\ 
& 3 & 0.21 & 0.22 & 0.21 & 0.21 &  --  \\ 
& 4 & 0.21 & 0.20 & 0.19 &  --  &  --  \\ 
& 5 & 0.29 & 0.29 & 0.29 &  --  &  --  \\ \midrule
\end{tabular}
\end{table}

\begin{table}[htb]
\caption{$h$-multigrid convergence factors for primal and flux coarsening applied to the LDG discretization of Poisson's equation on a uniform $n \times n \times n$ Cartesian grid as ${h\to 0}$ using MGPCG (see Figure 4.2).\vspace{-1.8em}}
\label{tab:3D_uniform_mgpcg}
\centering
\small
\begin{tabular}{@{}cl@{\qquad}ccccc@{}}
&&&&&& \\ \midrule
&& \multicolumn{5}{c}{$n$} \\
& $p$ & 4 & 8 & 16 & 32 & 64 \\ \cmidrule{1-7}
\multirow{5}{*}{Primal coarsening}
& 1 & 0.06 & 0.10 & 0.16 & 0.23 & 0.33 \\ 
& 2 & 0.07 & 0.10 & 0.14 & 0.20 & 0.28 \\ 
& 3 & 0.10 & 0.13 & 0.18 & 0.25 &  --  \\ 
& 4 & 0.09 & 0.12 & 0.17 &  --  &  --  \\ 
& 5 & 0.13 & 0.16 & 0.21 &  --  &  --  \\ \cmidrule{1-7}
\multirow{5}{*}{Flux coarsening}
& 1 & 0.06 & 0.08 & 0.09 & 0.09 & 0.10 \\ 
& 2 & 0.06 & 0.06 & 0.06 & 0.06 & 0.06 \\ 
& 3 & 0.09 & 0.09 & 0.09 & 0.09 &  --  \\ 
& 4 & 0.09 & 0.09 & 0.09 &  --  &  --  \\ 
& 5 & 0.12 & 0.12 & 0.12 &  --  &  --  \\ \midrule
\end{tabular}
\end{table}

\begin{table}[htb]
\caption{$p$-multigrid convergence factors for primal and flux coarsening applied to the LDG discretization of Poisson's equation on 2D uniform grids using MGPCG (see Figure 4.3).\vspace{-1.8em}}
\label{tab:2D_pmg}
\centering
\small
\begin{tabular}{@{}cl@{\qquad}cccc@{}}
&&&&& \\ \midrule
&& \multicolumn{4}{c}{$p$} \\
& $n$ & 1 & 2 & 4 & 8 \\ \cmidrule{1-6}
\multirow{8}{*}{Primal coarsening}
& 4 & 0.04 & 0.06 & 0.10 & 0.20 \\ 
& 8 & 0.09 & 0.12 & 0.19 & 0.36 \\ 
& 16 & 0.13 & 0.19 & 0.32 & 0.51 \\ 
& 32 & 0.20 & 0.30 & 0.45 & 0.63 \\ 
& 64 & 0.30 & 0.43 & 0.57 & 0.71 \\ 
& 128 & 0.40 & 0.52 & 0.67 & 0.79 \\ 
& 256 & 0.50 & 0.63 & 0.75 & 0.84 \\ 
& 512 & 0.60 & 0.71 & 0.80 &  --  \\ \cmidrule{1-6}
\multirow{8}{*}{Flux coarsening}
& 4 & 0.05 & 0.04 & 0.04 & 0.08 \\ 
& 8 & 0.06 & 0.05 & 0.05 & 0.08 \\ 
& 16 & 0.06 & 0.06 & 0.06 & 0.08 \\ 
& 32 & 0.07 & 0.06 & 0.07 & 0.08 \\ 
& 64 & 0.07 & 0.07 & 0.08 & 0.08 \\ 
& 128 & 0.07 & 0.08 & 0.08 & 0.08 \\ 
& 256 & 0.07 & 0.08 & 0.08 & 0.08 \\ 
& 512 & 0.08 & 0.08 & 0.08 &  --  \\ \midrule
\end{tabular}
\end{table}

\begin{table}[htb]
\caption{$p$-multigrid convergence factors for primal and flux coarsening applied to the LDG discretization of Poisson's equation on 3D uniform grids using MGPCG (see Figure 4.3).\vspace{-1.8em}}
\label{tab:3D_pmg}
\centering
\small
\begin{tabular}{@{}cl@{\qquad}cccc@{}}
&&&&& \\ \midrule
&& \multicolumn{4}{c}{$p$} \\
& $n$ & 1 & 2 & 4 & 8 \\ \cmidrule{1-6}
\multirow{5}{*}{Primal coarsening}
& 4 & 0.06 & 0.08 & 0.13 & 0.25 \\ 
& 8 & 0.10 & 0.15 & 0.26 & 0.41 \\ 
& 16 & 0.16 & 0.24 & 0.38 &  --  \\ 
& 32 & 0.23 & 0.35 & 0.50 &  --  \\ 
& 64 & 0.33 & 0.46 &  --  &  --  \\ \cmidrule{1-6}
\multirow{5}{*}{Flux coarsening}
& 4 & 0.06 & 0.04 & 0.06 & 0.12 \\ 
& 8 & 0.08 & 0.07 & 0.07 & 0.12 \\ 
& 16 & 0.09 & 0.08 & 0.08 &  --  \\ 
& 32 & 0.09 & 0.08 & 0.08 &  --  \\ 
& 64 & 0.10 & 0.08 &  --  &  --  \\ \midrule
\end{tabular}
\end{table}

\begin{table}[htb]
\caption{$h$-multigrid convergence factors for primal and flux coarsening applied to the LDG discretization of Poisson's equation on an adaptively refined grid in 2D using MGPCG (see Figure 4.6). The effective grid size $h_\text{eff}$ is the size of the smallest element on the finest grid.\vspace{-1.8em}}
\label{tab:2D_adaptive}
\centering
\small
\begin{tabular}{@{}cl@{\qquad}cccccc@{}}
&&&&&&& \\ \midrule
&& \multicolumn{6}{c}{$1/h_\text{eff}$} \\
& $p$ & 16 & 32 & 64 & 128 & 256 & 512 \\ \cmidrule{1-8}
\multirow{5}{*}{Primal coarsening}
& 1 & 0.03 & 0.09 & 0.17 & 0.25 & 0.37 & 0.47 \\ 
& 2 & 0.04 & 0.08 & 0.15 & 0.21 & 0.32 & 0.42 \\ 
& 3 & 0.07 & 0.11 & 0.18 & 0.25 & 0.35 & 0.44 \\ 
& 4 & 0.06 & 0.10 & 0.17 & 0.26 & 0.37 & 0.45 \\ 
& 5 & 0.09 & 0.13 & 0.19 & 0.29 & 0.38 & 0.47 \\ \cmidrule{1-8}
\multirow{5}{*}{Flux coarsening}
& 1 & 0.02 & 0.06 & 0.08 & 0.09 & 0.12 & 0.14 \\ 
& 2 & 0.03 & 0.05 & 0.09 & 0.11 & 0.13 & 0.14 \\ 
& 3 & 0.06 & 0.08 & 0.11 & 0.15 & 0.17 & 0.18 \\ 
& 4 & 0.06 & 0.08 & 0.13 & 0.18 & 0.19 & 0.21 \\ 
& 5 & 0.08 & 0.12 & 0.16 & 0.21 & 0.22 & 0.25 \\ \midrule
\end{tabular}
\end{table}

\begin{table}[htb]
\caption{$h$-multigrid convergence factors for primal and flux coarsening applied to the LDG discretization of Poisson's equation on an adaptively refined grid in 3D using MGPCG (see Figure 4.6). The effective grid size $h_\text{eff}$ is the size of the smallest element on the finest grid.\vspace{-1.8em}}
\label{tab:3D_adaptive}
\centering
\small
\begin{tabular}{@{}cl@{\qquad}ccccc@{}}
&&&&&& \\ \midrule
&& \multicolumn{5}{c}{$1/h_\text{eff}$} \\
& $p$ & 16 & 32 & 64 & 128 & 256 \\ \cmidrule{1-7}
\multirow{5}{*}{Primal coarsening}
& 1 & 0.05 & 0.11 & 0.19 & 0.28 & 0.39 \\ 
& 2 & 0.05 & 0.10 & 0.18 & 0.28 &  --  \\ 
& 3 & 0.08 & 0.14 & 0.21 & 0.30 &  --  \\ 
& 4 & 0.08 & 0.13 & 0.21 &  --  &  --  \\ 
& 5 & 0.11 & 0.16 &  --  &  --  &  --  \\ \cmidrule{1-7}
\multirow{5}{*}{Flux coarsening}
& 1 & 0.04 & 0.08 & 0.10 & 0.14 & 0.14 \\ 
& 2 & 0.04 & 0.07 & 0.14 & 0.18 &  --  \\ 
& 3 & 0.08 & 0.11 & 0.17 & 0.21 &  --  \\ 
& 4 & 0.08 & 0.11 & 0.21 &  --  &  --  \\ 
& 5 & 0.12 & 0.15 &  --  &  --  &  --  \\ \midrule
\end{tabular}
\end{table}

\begin{table}[htb]
\caption{$h$-multigrid convergence factors for flux coarsening and MGPCG applied to the LDG discretization of the single-phase Poisson problem on the curved domain shown in Figure~4.7 (see Figure 4.9). The grid size $n$ is the number of cells of the background Cartesian grid required to cover the longest extent of the domain.\vspace{-1.8em}}
\label{tab:lemniscate}
\centering
\small
\begin{tabular}{@{}cl@{\qquad}ccccccc@{}}
&&&&&&&& \\ \midrule
&& \multicolumn{7}{c}{$n$} \\
& $p$ & 16 & 32 & 64 & 128 & 256 & 512 & 1024 \\ \cmidrule{1-9}
& 1 & 0.01 & 0.05 & 0.06 & 0.08 & 0.09 & 0.10 & 0.11 \\ 
& 2 & 0.04 & 0.07 & 0.08 & 0.07 & 0.08 & 0.08 & 0.10 \\ 
& 3 & 0.06 & 0.10 & 0.11 & 0.11 & 0.12 & 0.11 & 0.11 \\ 
& 4 & 0.08 & 0.14 & 0.15 & 0.14 & 0.15 & 0.15 & 0.16 \\ 
& 5 & 0.10 & 0.16 & 0.19 & 0.18 & 0.21 & 0.18 & 0.18 \\ \midrule
\end{tabular}
\end{table}

\begin{table}[htb]
\caption{$h$-multigrid convergence factors for flux coarsening and MGPCG applied to the LDG discretization of the multi-phase elliptic interface test problem in (4.2) on the domain shown in Figure~4.8 (see Figure 4.9). The grid size $n$ is the number of cells of the background Cartesian grid required to cover the longest extent of the domain.\vspace{-1.8em}}
\label{tab:elliptic_interface}
\centering
\small
\begin{tabular}{@{}cl@{\qquad}ccccc@{}}
&&&&&& \\ \midrule
&& \multicolumn{5}{c}{$n$} \\
& $p$ & 64 & 128 & 256 & 512 & 1024 \\ \cmidrule{1-7}
& 1 &  --  &  --  & 0.17 & 0.20 & 0.27 \\ 
& 2 & 0.23 & 0.16 & 0.18 & 0.20 & 0.24 \\ 
& 3 & 0.19 & 0.17 & 0.17 & 0.21 & 0.28 \\ 
& 4 & 0.35 & 0.28 & 0.26 & 0.28 & 0.29 \\ \midrule
\end{tabular}
\end{table}